\documentclass{amsart}
\usepackage{amssymb}
\usepackage{amsfonts}

\newtheorem{theorem}{Theorem}
\theoremstyle{plain}

\newtheorem{definition}{Definition}

\newtheorem{lemma}{Lemma}

\newtheorem{problem}{Problem}

\newtheorem{remark}{Remark}

\numberwithin{equation}{section}
\input{tcilatex}

\begin{document}
\title[A geometric condition and necessity of energy]{A geometric condition,
necessity of energy, and two weight boundedness of fractional Riesz
transforms}
\author[E.T. Sawyer]{Eric T. Sawyer}
\address{ Department of Mathematics \& Statistics, McMaster University, 1280
Main Street West, Hamilton, Ontario, Canada L8S 4K1 }
\email{sawyer@mcmaster.ca}
\thanks{Research supported in part by NSERC}
\author[C.-Y. Shen]{Chun-Yen Shen}
\address{ Department of Mathematics \\
National Central University \\
Chungli, 32054, Taiwan }
\email{chunyshen@gmail.com}
\thanks{C.-Y. Shen supported in part by the NSC, through grant
NSC102-2115-M-008-015-MY2}
\author[I. Uriarte-Tuero]{Ignacio Uriarte-Tuero}
\address{ Department of Mathematics \\
Michigan State University \\
East Lansing MI }
\email{ignacio@math.msu.edu}
\thanks{ I. Uriarte-Tuero has been partially supported by grants DMS-1056965
(US NSF), MTM2010-16232, MTM2009-14694-C02-01 (Spain), and a Sloan
Foundation Fellowship. }
\date{March 16, 2014}

\begin{abstract}
Let $\sigma $ and $\omega $ be locally finite positive Borel measures on $%
\mathbb{R}^{n}$ with no common point masses. We assume that at least one of
the two measures $\sigma $ and $\omega $ is supported on a line in $\mathbb{R%
}^{n}$. Let $\mathbf{R}^{\alpha ,n}$\ be the $\alpha $-fractional Riesz
transform vector on $\mathbb{R}^{n}$. We prove that the energy conditions in 
\textit{arXiv:1302.5093v7} are implied by the $\mathcal{A}_{2}^{\alpha }$
and cube testing conditions for $\mathbf{R}^{\alpha ,n}$. Then we apply the
main theorem there to give a T1 theorem for $\mathbf{R}^{\alpha ,n}$: namely
that $\mathbf{R}^{\alpha ,n}$ is bounded from $L^{2}\left( \sigma \right) $
to $L^{2}\left( \omega \right) $ \emph{if and only if} the $\mathcal{A}%
_{2}^{\alpha }$ conditions hold, the cube testing conditions for $\mathbf{R}%
^{\alpha ,n}$\textbf{\ }and its dual both hold, and the weak boundedness
property for $\mathbf{R}^{\alpha ,n}$\ holds.
\end{abstract}

\maketitle
\tableofcontents

\section{Introduction}

In \cite{SaShUr}, under a side assumption that certain \emph{energy
conditions} hold,\ the authors show in particular that the two weight
inequality%
\begin{equation}
\left\Vert \mathbf{R}^{\alpha ,n}\left( f\sigma \right) \right\Vert
_{L^{2}\left( \omega \right) }\lesssim \left\Vert f\right\Vert _{L^{2}\left(
\sigma \right) },  \label{2 weight}
\end{equation}%
for the vector of Riesz transforms $\mathbf{R}^{\alpha ,n}$ in $\mathbb{R}%
^{n}$ (with $0\leq \alpha <n$) holds if and only if the $\mathcal{A}_{2}$
conditions hold, the cube testing conditions hold, and the weak boundedness
property holds. It is not known at the time of this writing whether or not
these or any other energy conditions are necessary for \emph{any} vector $%
\mathbf{T}^{\alpha ,n}$ of fractional singular integrals in $\mathbb{R}^{n}$
with $n\geq 2$, apart from the trivial case of positive operators. In
particular there are no known counterexamples. We also showed in \cite%
{SaShUr2} and \cite{SaShUr3} that the technique of reversing energy,
typically used to prove energy conditions, fails spectacularly in higher
dimension (and we thank M. Lacey for showing us this failure for the Cauchy
transform with the circle measure). See also the counterexamples for the
fractional Riesz transforms in \cite{LaWi2}.

The purpose of this paper is to show that if $\sigma $ and $\omega $ are
locally finite positive Borel measures without common point masses, and at
least \emph{one} of the two measures $\sigma $ and $\omega $ is supported on
a line in $\mathbb{R}^{n}$, then the energy conditions are indeed necessary
for boundedness of the fractional Riesz transform $\mathbf{R}^{\alpha ,n}$,
and hence that a T1 theorem holds for $\mathbf{R}^{\alpha ,n}$. M. Lacey and
B. Wick \cite{LaWi} have independently obtained a similar result for the
Cauchy transform in the plane, and the five authors have combined on the
paper \cite{LaSaShUrWi}. The vector of $\alpha $-fractional Riesz transforms
is given by%
\begin{equation*}
\mathbf{R}^{\alpha ,n}=\left\{ R_{\ell }^{\alpha ,n}:1\leq \ell \leq
n\right\} ,\ \ \ \ \ 0\leq \alpha <n,
\end{equation*}%
where the component Riesz transforms $R_{\ell }^{\alpha ,n}$ are the
convolution fractional singular integrals $R_{\ell }^{\alpha ,n}f\equiv
K_{\ell }^{\alpha ,n}\ast f$ with odd kernel defined by%
\begin{equation*}
K_{\ell }^{\alpha ,n}\left( w\right) \equiv c_{\alpha ,n}\frac{w^{\ell }}{%
\left\vert w\right\vert ^{n+1-\alpha }}.
\end{equation*}%
Finally, we remark that the T1 theorem under this geometric condition has
application to the weighted discrete Hilbert transform $H_{\left( \Gamma
,v\right) }$ when the sequence $\Gamma $ is supported on a line in the
complex plane. See \cite{BeMeSe} where $H_{\left( \Gamma ,v\right) }$ is
essentially the Cauchy transform with $n=2$ and $\alpha =1$.

We now recall a special case of our main two weight theorem from \cite%
{SaShUr}. Let $\mathcal{Q}^{n}$ denote the collection of all cubes in $%
\mathbb{R}^{n}$, and denote by $\mathcal{D}^{n}$ a dyadic grid in $\mathbb{R}%
^{n}$. The definitions of the remaining terms used below will be given in
the next section.

\begin{theorem}
\label{T1 theorem}Suppose that $\mathbf{R}^{\alpha ,n}$ is the vector of $%
\alpha $-fractional Riesz transforms in $\mathbb{R}^{n}$, and that $\omega $
and $\sigma $ are positive Borel measures on $\mathbb{R}^{n}$ without common
point masses. Set $\mathbf{R}_{\sigma }^{\alpha ,n}f=\mathbf{R}^{\alpha
,n}\left( f\sigma \right) $ for any smooth truncation of $\mathbf{R}^{\alpha
,n}$.

\begin{enumerate}
\item Suppose $0\leq \alpha <n$ and that $\gamma \geq 2$ is given. Then the
operator $\mathbf{R}_{\sigma }^{\alpha ,n}$ is bounded from $L^{2}\left(
\sigma \right) $ to $L^{2}\left( \omega \right) $, i.e. 
\begin{equation}
\left\Vert \mathbf{R}_{\sigma }^{\alpha ,n}f\right\Vert _{L^{2}\left( \omega
\right) }\leq \mathfrak{N}_{\mathbf{R}_{\sigma }^{\alpha ,n}}\left\Vert
f\right\Vert _{L^{2}\left( \sigma \right) },  \label{two weight}
\end{equation}%
uniformly in smooth truncations of $T^{\alpha }$, and moreover%
\begin{equation*}
\mathfrak{N}_{\mathbf{R}_{\sigma }^{\alpha ,n}}\leq C_{\alpha }\left( \sqrt{%
\mathcal{A}_{2}^{\alpha }+\mathcal{A}_{2}^{\alpha ,\ast }}+\mathfrak{T}_{%
\mathbf{R}_{\sigma }^{\alpha ,n}}+\mathfrak{T}_{\mathbf{R}_{\sigma }^{\alpha
,n}}^{\ast }+\mathcal{E}_{\alpha }+\mathcal{E}_{\alpha }^{\ast }+\mathcal{WBP%
}_{\mathbf{R}_{\sigma }^{\alpha ,n}}\right) ,
\end{equation*}%
provided that the two dual $\mathcal{A}_{2}^{\alpha }$ conditions hold, and
the two dual testing conditions for $\mathbf{R}_{\sigma }^{\alpha ,n}$ hold,
the weak boundedness property for $\mathbf{R}_{\sigma }^{\alpha ,n}$ holds
for a sufficiently large constant $C$ depending on the goodness parameter $%
\mathbf{r}$, and provided that the two dual energy conditions $\mathcal{E}%
_{\alpha }+\mathcal{E}_{\alpha }^{\ast }<\infty $ hold uniformly over all
dyadic grids $\mathcal{D}^{n}$, and where the goodness parameters $\mathbf{r}
$ and $\varepsilon $ implicit in the definition of $\mathcal{M}_{\mathbf{r}-%
\limfunc{deep}}^{\ell }\left( K\right) $ are fixed sufficiently large and
small respectively depending on $n$, $\alpha $ and $\gamma $.

\item Conversely, suppose $0\leq \alpha <n$ and that the Riesz transform
vector $\mathbf{R}_{\sigma }^{\alpha ,n}$ is bounded from $L^{2}\left(
\sigma \right) $ to $L^{2}\left( \omega \right) $, 
\begin{equation*}
\left\Vert \mathbf{R}_{\sigma }^{\alpha ,n}f\right\Vert _{L^{2}\left( \omega
\right) }\leq \mathfrak{N}_{\mathbf{R}_{\sigma }^{\alpha ,n}}\left\Vert
f\right\Vert _{L^{2}\left( \sigma \right) }.
\end{equation*}%
Then the testing conditions and weak boundedness property hold for $\mathbf{R%
}_{\sigma }^{\alpha ,n}$, the fractional $\mathcal{A}_{2}^{\alpha }$
conditions hold, and moreover,%
\begin{equation*}
\sqrt{\mathcal{A}_{2}^{\alpha }+\mathcal{A}_{2}^{\alpha ,\ast }}+\mathfrak{T}%
_{\mathbf{R}_{\sigma }^{\alpha ,n}}+\mathfrak{T}_{\mathbf{R}_{\sigma
}^{\alpha ,n}}^{\ast }+\mathcal{WBP}_{\mathbf{R}_{\sigma }^{\alpha ,n}}\leq C%
\mathfrak{N}_{\mathbf{R}_{\sigma }^{\alpha ,n}}.
\end{equation*}
\end{enumerate}
\end{theorem}

\begin{problem}
It is an open question whether or not the energy conditions are necessary
for boundedness of $\mathbf{R}_{\sigma }^{\alpha ,n}$. See \cite{SaShUr3}
for a failure of \emph{energy reversal} in higher dimensions - such an
energy reversal was used in dimension $n=1$ to prove the necessity of the
energy condition for the Hilbert transform.
\end{problem}

\begin{remark}
The boundedness of an individual operator $T^{\alpha }$ cannot in general
imply the finiteness of either $A_{2}^{\alpha }$ or $\mathcal{E}_{\alpha }$.
For a trivial example, if $\sigma $ and $\omega $ are supported on the $x$%
-axis in the plane, then the second Riesz tranform $R_{2}$ is the zero
operator from $L^{2}\left( \sigma \right) $ to $L^{2}\left( \omega \right) $%
, simply because the kernel $K_{2}\left( x,y\right) $ of $R_{2}$ satisfies $%
K_{2}\left( \left( x_{1},0\right) ,\left( y_{1},0\right) \right) =\frac{0-0}{%
\left\vert x_{1}-y_{1}\right\vert ^{3-\alpha }}=0$.
\end{remark}

\begin{remark}
\label{surgery}In \cite{LaWi2}, M. Lacey and B. Wick use the NTV technique
of surgery to show that the weak boundedness property for the Riesz
transform vector $\mathbf{R}^{\alpha ,n}$ is implied by the $\mathcal{A}%
_{2}^{\alpha }$ and cube testing conditions, and this has the consequence of
eliminating the weak boundedness property as a condition from the statement
of Theorem \ref{T1 theorem}.
\end{remark}

The next result shows that the energy conditions are in fact necessary for
boundedness of the Riesz transform vector $\mathbf{R}_{\sigma }^{\alpha ,n}$
when one of the measures is supported on a line.

\begin{theorem}
\label{main}Let $\sigma $ and $\omega $ be locally finite positive Borel
measures on $\mathbb{R}^{n}$ with no common point masses. Suppose that $%
\mathbf{R}^{\alpha ,n}$ is the fractional Riesz transform with $0\leq \alpha
<n$, and consider the tangent line truncations for $\mathbf{R}^{\alpha ,n}$
in the testing conditions. If at least one of the measures $\sigma $ and $%
\omega $ is supported on a line, then%
\begin{equation*}
\mathcal{E}_{\alpha }\lesssim \sqrt{\mathcal{A}_{2}^{\alpha }}+\mathfrak{T}_{%
\mathbf{R}^{\alpha ,n}}\text{ and }\mathcal{E}_{\alpha }^{\ast }\lesssim 
\sqrt{\mathcal{A}_{2}^{\alpha ,\ast }}+\mathfrak{T}_{\mathbf{R}^{\alpha
,n}}^{\ast }.
\end{equation*}
\end{theorem}

If we combine Theorems \ref{main} and \ref{T1 theorem}, we obtain the
following theorem as a corollary, which generalizes the T1 theorem for the
Hilbert transform (\cite{Lac}, \cite{LaSaShUr3}). See also related work in
the references given at the end of the paper. We use notation as in Theorem %
\ref{T1 theorem}.

\begin{theorem}
\label{final}Let $\sigma $ and $\omega $ be locally finite positive Borel
measures on $\mathbb{R}^{n}$ with no common point masses. Suppose that $%
\mathbf{R}^{\alpha ,n}$ is the fractional Riesz transform with $0\leq \alpha
<n$. Set $\mathbf{R}_{\sigma }^{\alpha ,n}f=\mathbf{R}^{\alpha ,n}\left(
f\sigma \right) $ for any smooth truncation of $\mathbf{R}^{\alpha ,n}$. If
at least one of the measures $\sigma $ and $\omega $ is supported on a line,
then the operator norm $\mathfrak{N}_{\mathbf{R}^{\alpha ,n}}$ of $\mathbf{R}%
_{\sigma }^{\alpha ,n}$ as an operator from\thinspace $L^{2}\left( \sigma
\right) $ to $L^{2}\left( \omega \right) $, uniformly in smooth truncations,
satisfies%
\begin{equation*}
\mathfrak{N}_{\mathbf{R}^{\alpha ,n}}\approx C_{\alpha }\left( \sqrt{%
\mathcal{A}_{2}^{\alpha }+\mathcal{A}_{2}^{\alpha ,\ast }}+\mathfrak{T}_{%
\mathbf{R}^{\alpha ,n}}+\mathfrak{T}_{\mathbf{R}^{\alpha ,n}}^{\ast }+%
\mathcal{WBP}_{\mathbf{R}^{\alpha ,n}}\right) .
\end{equation*}
\end{theorem}

\section{Definitions}

As mentioned above, the $\alpha $-fractional \emph{Riesz} vector $\mathbf{R}%
^{\alpha ,n}=\left\{ R_{\ell }^{\alpha ,n}:1\leq \ell \leq n\right\} $ has
as components the Riesz transforms $R_{\ell }^{n,\alpha }$ with odd kernel $%
K_{\ell }^{\alpha ,n}\left( w\right) =\frac{\Omega _{\ell }\left( w\right) }{%
\left\vert w\right\vert ^{n-\alpha }}$. The \emph{tangent line truncation} of%
\emph{\ }the Riesz transform $R_{\ell }^{\alpha ,n}$ has kernel $\Omega
_{\ell }\left( w\right) \psi _{\delta ,R}^{\alpha }\left( \left\vert
w\right\vert \right) $ where $\psi _{\delta ,R}^{\alpha }$ is continuously
differentiable on an interval $\left( 0,S\right) $ with $0<\delta <R<S$, and
where $\psi _{\delta ,R}^{\alpha }\left( r\right) =r^{\alpha -n}$ if $\delta
\leq r\leq R$, and has constant derivative on both $\left( 0,\delta \right) $
and $\left( R,S\right) $ where $\psi _{\delta ,R}^{\alpha }\left( S\right) =0
$. As shown in the one dimensional case in \cite{LaSaShUr3}, boundedness of $%
R_{\ell }^{n,\alpha }$ with one set of appropriate truncations together with
the $\mathcal{A}_{2}^{\alpha }$ condition below, is equivalent to
boundedness of $R_{\ell }^{n,\alpha }$ with all truncations.

\subsection{Cube testing, the weak boundedness property, and the $\mathcal{A}%
_{2}^{\protect\alpha }$ conditions}

The following `dual' cube testing conditions are necessary for the
boundedness of $\mathbf{R}^{\alpha ,n}$ from $L^{2}\left( \sigma \right) $
to $L^{2}\left( \omega \right) $:%
\begin{eqnarray*}
\mathfrak{T}_{\mathbf{R}^{\alpha ,n}}^{2} &\equiv &\sup_{Q\in \mathcal{Q}%
^{n}}\frac{1}{\left\vert Q\right\vert _{\sigma }}\int_{Q}\left\vert \mathbf{R%
}^{\alpha ,n}\left( \mathbf{1}_{Q}\sigma \right) \right\vert ^{2}\omega
<\infty , \\
\left( \mathfrak{T}_{\mathbf{R}^{\alpha ,n}}^{\ast }\right) ^{2} &\equiv
&\sup_{Q\in \mathcal{Q}^{n}}\frac{1}{\left\vert Q\right\vert _{\omega }}%
\int_{Q}\left\vert \left( \mathbf{R}^{\alpha ,n}\right) ^{\ast }\left( 
\mathbf{1}_{Q}\omega \right) \right\vert ^{2}\sigma <\infty .
\end{eqnarray*}%
The weak boundedness property for $\mathbf{R}^{\alpha ,n}$ with constant $C$
is given by 
\begin{eqnarray*}
&&\left\vert \int_{Q}\mathbf{R}^{\alpha ,n}\left( 1_{Q^{\prime }}\sigma
\right) d\omega \right\vert \leq \mathcal{WBP}_{\mathbf{R}^{\alpha ,n}}\sqrt{%
\left\vert Q\right\vert _{\omega }\left\vert Q^{\prime }\right\vert _{\sigma
}}, \\
&&\ \ \ \ \ \text{for all cubes }Q,Q^{\prime }\text{ with }\frac{1}{C}\leq 
\frac{\left\vert Q\right\vert ^{\frac{1}{n}}}{\left\vert Q^{\prime
}\right\vert ^{\frac{1}{n}}}\leq C, \\
&&\ \ \ \ \ \text{and either }Q\subset 3Q^{\prime }\setminus Q^{\prime }%
\text{ or }Q^{\prime }\subset 3Q\setminus Q.
\end{eqnarray*}%
Now let $\mu $ be a locally finite positive Borel measure on $\mathbb{R}^{n}$%
, and suppose $Q$ is a cube in $\mathbb{R}^{n}$. The two $\alpha $%
-fractional Poisson integrals of $\mu $ on a cube $Q$ are given by:%
\begin{eqnarray*}
\mathrm{P}^{\alpha }\left( Q,\mu \right) &\equiv &\int_{\mathbb{R}^{n}}\frac{%
\left\vert Q\right\vert ^{\frac{1}{n}}}{\left( \left\vert Q\right\vert ^{%
\frac{1}{n}}+\left\vert x-x_{Q}\right\vert \right) ^{n+1-\alpha }}d\mu
\left( x\right) , \\
\mathcal{P}^{\alpha }\left( Q,\mu \right) &\equiv &\int_{\mathbb{R}%
^{n}}\left( \frac{\left\vert Q\right\vert ^{\frac{1}{n}}}{\left( \left\vert
Q\right\vert ^{\frac{1}{n}}+\left\vert x-x_{Q}\right\vert \right) ^{2}}%
\right) ^{n-\alpha }d\mu \left( x\right) .
\end{eqnarray*}%
We refer to $\mathrm{P}^{\alpha }$ as the \emph{standard} Poisson integral
and to $\mathcal{P}^{\alpha }$ as the \emph{reproducing} Poisson integral.
Let $\sigma $ and $\omega $ be locally finite positive Borel measures on $%
\mathbb{R}^{n}$ with no common point masses, and suppose $0\leq \alpha <n$.
The classical $A_{2}^{\alpha }$ constant is defined by 
\begin{equation*}
A_{2}^{\alpha }\equiv \sup_{Q\in \mathcal{Q}^{n}}\frac{\left\vert
Q\right\vert _{\sigma }}{\left\vert Q\right\vert ^{1-\frac{\alpha }{n}}}%
\frac{\left\vert Q\right\vert _{\omega }}{\left\vert Q\right\vert ^{1-\frac{%
\alpha }{n}}},
\end{equation*}%
and the one-sided constants $\mathcal{A}_{2}^{\alpha }$ and $\mathcal{A}%
_{2}^{\alpha ,\ast }$ for the weight pair $\left( \sigma ,\omega \right) $
are defined by%
\begin{eqnarray*}
\mathcal{A}_{2}^{\alpha } &\equiv &\sup_{Q\in \mathcal{Q}^{n}}\mathcal{P}%
^{\alpha }\left( Q,\sigma \right) \frac{\left\vert Q\right\vert _{\omega }}{%
\left\vert Q\right\vert ^{1-\frac{\alpha }{n}}}<\infty , \\
\mathcal{A}_{2}^{\alpha ,\ast } &\equiv &\sup_{Q\in \mathcal{Q}^{n}}\mathcal{%
P}^{\alpha }\left( Q,\omega \right) \frac{\left\vert Q\right\vert _{\sigma }%
}{\left\vert Q\right\vert ^{1-\frac{\alpha }{n}}}<\infty .
\end{eqnarray*}

\subsection{Energy conditions}

We begin by briefly recalling some of the notation used in \cite{SaShUr}.
Given a dyadic cube $K\in \mathcal{D}$ and a positive measure $\mu $ we
define the Haar projection $\mathsf{P}_{K}^{\mu }\equiv \sum_{_{J\in 
\mathcal{D}:\ J\subset K}}\bigtriangleup _{J}^{\mu }$. Now we recall the
definition of a \emph{good} dyadic cube - see \cite{NTV3} and \cite{LaSaUr2}
for more detail.

\begin{definition}
Let $\mathbf{r}\in \mathbb{N}$ and $0<\varepsilon <1$. A dyadic cube $J$ is $%
\left( \mathbf{r},\varepsilon \right) $\emph{-good}, or simply \emph{good},
if for \emph{every} dyadic supercube $I$, it is the case that \textbf{either}
$J$ has side length at least $2^{-\mathbf{r}}$ times that of $I$, \textbf{or}
$J\Subset _{\mathbf{r}}I$ is $\left( \mathbf{r},\varepsilon \right) $-deeply
embedded in $I$.
\end{definition}

Here we say that a dyadic cube $J$ is $\left( \mathbf{r},\varepsilon \right) 
$-\emph{deeply embedded} in a dyadic cube $K$, or simply $\mathbf{r}$\emph{%
-deeply embedded} in $K$, which we write as $J\Subset _{\mathbf{r}}K$, when $%
J\subset K$ and both 
\begin{eqnarray}
\left\vert J\right\vert ^{\frac{1}{n}} &\leq &2^{-\mathbf{r}}\left\vert
K\right\vert ^{\frac{1}{n}},  \label{def deep embed} \\
\limfunc{dist}\left( J,\partial K\right)  &\geq &\frac{1}{2}\left\vert
J\right\vert ^{\frac{\varepsilon }{n}}\left\vert K\right\vert ^{\frac{%
1-\varepsilon }{n}}.  \notag
\end{eqnarray}%
We say that $J$ is $\mathbf{r}$\emph{-nearby} in $K$ when $J\subset K$ and%
\begin{equation*}
\left\vert J\right\vert ^{\frac{1}{n}}>2^{-\mathbf{r}}\left\vert
K\right\vert ^{\frac{1}{n}}.
\end{equation*}%
We denote the set of such good dyadic cubes by $\mathcal{D}_{\limfunc{good}}$%
.

Then we define the smaller `good' Haar projection $\mathsf{P}_{K}^{\limfunc{%
good},\omega }$ by%
\begin{equation*}
\mathsf{P}_{K}^{\limfunc{good},\mu }f\equiv \sum_{_{J\in \mathcal{G}\left(
K\right) }}\bigtriangleup _{J}^{\mu }f,
\end{equation*}%
where $\mathcal{G}\left( K\right) $ consists of the good subcubes of $K$:%
\begin{equation*}
\mathcal{G}\left( K\right) \equiv \left\{ J\in \mathcal{D}_{\limfunc{good}%
}:J\subset K\right\} ,
\end{equation*}%
and also the larger `subgood' Haar projection $\mathsf{P}_{K}^{\limfunc{%
subgood},\mu }$ by%
\begin{equation*}
\mathsf{P}_{K}^{\limfunc{subgood},\mu }f\equiv \sum_{_{J\in \mathcal{M}_{%
\limfunc{good}}\left( K\right) }}\sum_{J^{\prime }\subset J}\bigtriangleup
_{J^{\prime }}^{\mu }f,
\end{equation*}%
where $\mathcal{M}_{\limfunc{good}}\left( K\right) $ consists of the \emph{%
maximal} good subcubes of $K$. We thus have 
\begin{eqnarray*}
\left\Vert \mathsf{P}_{K}^{\limfunc{good},\mu }\mathbf{x}\right\Vert
_{L^{2}\left( \mu \right) }^{2} &\leq &\left\Vert \mathsf{P}_{K}^{\limfunc{%
subgood},\mu }\mathbf{x}\right\Vert _{L^{2}\left( \mu \right) }^{2} \\
&\leq &\left\Vert \mathsf{P}_{I}^{\mu }\mathbf{x}\right\Vert _{L^{2}\left(
\mu \right) }^{2}=\int_{I}\left\vert \mathbf{x}-\left( \frac{1}{\left\vert
I\right\vert _{\mu }}\int_{I}\mathbf{x}dx\right) \right\vert ^{2}d\mu \left(
x\right) ,\ \ \ \ \ \mathbf{x}=\left( x_{1},...,x_{n}\right) ,
\end{eqnarray*}%
where $\mathsf{P}_{I}^{\mu }\mathbf{x}$ is the orthogonal projection of the
identity function $\mathbf{x}:\mathbb{R}^{n}\rightarrow \mathbb{R}^{n}$ onto
the vector-valued subspace of $\oplus _{k=1}^{n}L^{2}\left( \mu \right) $
consisting of functions supported in $I$ with $\mu $-mean value zero.

We use the collection $\mathcal{M}_{\mathbf{r}-\limfunc{deep}}\left(
K\right) $ of \emph{maximal} $\mathbf{r}$-deeply embedded dyadic subcubes of
a dyadic cube $K$. We let $J^{\ast }=\gamma J$ where $\gamma \geq 2$. The
goodness parameter $\mathbf{r}$ is chosen sufficiently large, depending on $%
\varepsilon $ and $\gamma $, that the bounded overlap property 
\begin{equation}
\sum_{J\in \mathcal{M}_{\mathbf{r}-\limfunc{deep}}\left( K\right) }\mathbf{1}%
_{J^{\ast }}\leq \beta \mathbf{1}_{K}\ ,  \label{bounded overlap}
\end{equation}%
holds for some positive constant $\beta $ depending only on $n,\gamma ,%
\mathbf{r}$ and $\varepsilon $. We will also need the following refinement
of $\mathcal{M}_{\mathbf{r}-\limfunc{deep}}\left( K\right) $ for each $\ell
\geq 0$ that consists of some of the maximal cubes $Q$, whose $\ell $-fold
dyadic parent $\pi ^{\ell }Q$ is $\mathbf{r}$-deeply embedded in $K$:%
\begin{equation*}
\mathcal{M}_{\mathbf{r}-\limfunc{deep}}^{\ell }\left( K\right) \equiv
\left\{ J\in \mathcal{M}_{\mathbf{r}-\limfunc{deep}}\left( \pi ^{\ell
}K\right) :J\subset L\text{ for some }L\in \mathcal{M}_{\limfunc{deep}%
}\left( K\right) \right\} .
\end{equation*}%
Since $J\in \mathcal{M}_{\mathbf{r}-\limfunc{deep}}^{\ell }\left( K\right) $
implies $\gamma J\subset K$, we also have from (\ref{bounded overlap}) that%
\begin{equation}
\sum_{J\in \mathcal{M}_{\mathbf{r}-\limfunc{deep}}^{\left( \ell \right)
}\left( K\right) }\mathbf{1}_{J^{\ast }}\leq \beta \mathbf{1}_{K}\ ,\ \ \ \
\ \text{for each }\ell \geq 0.  \label{bounded overlap'}
\end{equation}%
Of course $\mathcal{M}_{\mathbf{r}-\limfunc{deep}}^{0}\left( K\right) =%
\mathcal{M}_{\mathbf{r}-\limfunc{deep}}\left( K\right) $, but $\mathcal{M}_{%
\mathbf{r}-\limfunc{deep}}^{\ell }\left( K\right) $ is in general a finer
subdecomposition of $K$ the larger $\ell $ is, and may in fact be empty. The
following definition of the energy constant $\mathcal{E}_{\alpha }$ is
larger than that used in \cite{SaShUr}.

\begin{definition}
\label{energy condition}Suppose $\sigma $ and $\omega $ are positive Borel
measures on $\mathbb{R}^{n}$ without common point masses. Then the energy
condition constant $\mathcal{E}_{\alpha }$ is given by%
\begin{equation*}
\left( \mathcal{E}_{\alpha }\right) ^{2}\equiv \sup_{\ell \geq 0}\sup_{I=%
\dot{\cup}I_{r}}\frac{1}{\left\vert I\right\vert _{\sigma }}%
\sum_{r=1}^{\infty }\sum_{J\in \mathcal{M}_{\mathbf{r}-\limfunc{deep}}^{\ell
}\left( I_{r}\right) }\left( \frac{\mathrm{P}^{\alpha }\left( J,\mathbf{1}%
_{I\setminus \gamma J}\sigma \right) }{\left\vert J\right\vert ^{\frac{1}{n}}%
}\right) ^{2}\left\Vert \mathsf{P}_{J}^{\limfunc{subgood},\omega }\mathbf{x}%
\right\Vert _{L^{2}\left( \omega \right) }^{2}\ ,
\end{equation*}%
where $\sup_{I=\dot{\cup}I_{r}}$ above is taken over

\begin{enumerate}
\item all dyadic grids $\mathcal{D}$,

\item all $\mathcal{D}$-dyadic cubes $I$,

\item and all subpartitions $\left\{ I_{r}\right\} _{r=1}^{\infty }$ of the
cube $I$ into $\mathcal{D}$-dyadic subcubes $I_{r}$.
\end{enumerate}
\end{definition}

There is a similar definition for the dual (backward) energy condition that
simply interchanges $\sigma $ and $\omega $ everywhere. These definitions
of\ the energy conditions depend on the choice of goodness parameters $%
\mathbf{r}$ and $\varepsilon $. We can `plug the $\gamma $-hole' in the
Poisson integral $\mathrm{P}^{\alpha }\left( J,\mathbf{1}_{I\setminus \gamma
J}\sigma \right) $ using the $A_{2}^{\alpha }$ condition and the bounded
overlap property (\ref{bounded overlap'}). Indeed, with 
\begin{equation}
\left( \mathcal{E}_{\alpha }^{\limfunc{plug}}\right) ^{2}\equiv \sup_{\ell
\geq 0}\sup_{I=\dot{\cup}I_{r}}\frac{1}{\left\vert I\right\vert _{\sigma }}%
\sum_{r=1}^{\infty }\sum_{J\in \mathcal{M}_{\mathbf{r}-\limfunc{deep}}^{\ell
}\left( I_{r}\right) }\left( \frac{\mathrm{P}^{\alpha }\left( J,\mathbf{1}%
_{I}\sigma \right) }{\left\vert J\right\vert ^{\frac{1}{n}}}\right)
^{2}\left\Vert \mathsf{P}_{J}^{\limfunc{subgood},\omega }\mathbf{x}%
\right\Vert _{L^{2}\left( \omega \right) }^{2}\ ,  \label{plug}
\end{equation}%
we have, as shown in \cite{SaShUr}, that 
\begin{equation}
\left( \mathcal{E}_{\alpha }^{\limfunc{plug}}\right) ^{2}\lesssim \left( 
\mathcal{E}_{\alpha }\right) ^{2}+\beta A_{2}^{\alpha },
\label{plug the hole refined}
\end{equation}%
upon using (\ref{bounded overlap'}).

\subsection{Energy lemma}

We will need the following elementary special case of the Energy Lemma from 
\cite{SaShUr}.

\begin{lemma}[\textbf{Energy Lemma}]
\label{ener}Let $J\ $be a cube in $\mathcal{D}^{\omega }$. Let $\Psi _{J}$
be an $L^{2}\left( \omega \right) $ function supported in $J$ and with $%
\omega $-integral zero. Let $\nu $ be a positive measure supported in $%
\mathbb{R}^{n}\setminus \gamma J$ with $\gamma \geq 2$. Then we have%
\begin{equation*}
\left\vert \left\langle \mathbf{R}^{\alpha ,n}\left( \nu \right) ,\Psi
_{J}\right\rangle _{\omega }\right\vert \lesssim \left\Vert \Psi
_{J}\right\Vert _{L^{2}\left( \omega \right) }\left( \frac{\mathrm{P}%
^{\alpha }\left( J,\nu \right) }{\left\vert J\right\vert ^{\frac{1}{n}}}%
\right) \left\Vert \mathsf{P}_{J}^{\omega }\mathbf{x}\right\Vert
_{L^{2}\left( \omega \right) }.
\end{equation*}
\end{lemma}

\section{One measure supported in a line}

In this section we prove Theorem \ref{main}, i.e. we prove the necessity of
the energy conditions for the $\mathcal{A}_{2}^{\alpha }$ conditions and the
testing conditions $\mathfrak{T}_{\mathbf{R}^{\alpha ,n}}$ and $\mathfrak{T}%
_{\mathbf{R}^{\alpha ,n}}^{\ast }$ associated to the tangent line
truncations of the $\alpha $-fractional Riesz transform $\mathbf{R}^{\alpha
,n}$, when just \emph{one} of the measures $\sigma $ or $\omega $ is
supported in a line $L$, and the other measure is arbitrary. The
one-dimensional character of just one of the measures is enough to
circumvent the failure of strong reversal of energy as described in \cite%
{SaShUr2} and \cite{SaShUr3}.

Fix a dyadic grid $\mathcal{D}$, and suppose that $\omega $ is supported in
a line $L$. We will show that both energy conditions hold relative to $%
\mathcal{D}$. We can suppose that $L$ is the $x_{1}$-axis, since using that
the Riesz transform vector $\mathbf{R}^{\alpha ,n}$ is rotation invariant,
one can verify that the argument below does not depend in a critical way on
this or any other special relationship between $\mathcal{D}$ and $L$.

\subsection{Backward energy condition}

The dual (backward) energy condition $\mathcal{E}_{\alpha }^{\ast }\lesssim 
\mathfrak{T}_{\mathbf{R}^{\alpha ,n}}^{\ast }+\sqrt{\mathcal{A}_{2}^{\alpha
,\ast }}$ is the more straightforward of the two to verify, and so we turn
to it first. We must show%
\begin{equation*}
\sup_{\ell \geq 0}\sum_{r=1}^{\infty }\sum_{J\in \mathcal{M}_{\limfunc{deep}%
}^{\ell }\left( I_{r}\right) }\left( \frac{\mathrm{P}^{\alpha }\left( J,%
\mathbf{1}_{I\setminus J^{\ast }}\omega \right) }{\left\vert J\right\vert ^{%
\frac{1}{n}}}\right) ^{2}\left\Vert \mathsf{P}_{J}^{\limfunc{subgood},\sigma
}\mathbf{x}\right\Vert _{L^{2}\left( \sigma \right) }^{2}\leq \left( \left( 
\mathfrak{T}_{\mathbf{R}^{\alpha ,n}}^{\ast }\right) ^{2}+\mathcal{A}%
_{2}^{\alpha ,\ast }\right) \left\vert I\right\vert _{\omega }\ ,
\end{equation*}%
for all partitions of a dyadic cube $I=\dbigcup\limits_{r=1}^{\infty }I_{r}$
into dyadic subcubes $I_{r}$. We fix $\ell \geq 0$ and suppress both $\ell $
and $\mathbf{r}$ in the notation $\mathcal{M}_{\limfunc{deep}}\left(
I_{r}\right) =\mathcal{M}_{\mathbf{r}-\limfunc{deep}}^{\ell }\left(
I_{r}\right) $. Recall that $J^{\ast }=\gamma J$, and that the bounded
overlap property (\ref{bounded overlap'}) holds. We may of course assume
that $I$ intersects the $x_{1}$-axis $L$. Now we set $\mathcal{M}_{\limfunc{%
deep}}\equiv \dbigcup\limits_{r=1}^{\infty }\mathcal{M}_{\limfunc{deep}%
}\left( I_{r}\right) $ and write%
\begin{equation*}
\sum_{r=1}^{\infty }\sum_{J\in \mathcal{M}_{\limfunc{deep}}\left(
I_{r}\right) }\left( \frac{\mathrm{P}^{\alpha }\left( J,\mathbf{1}%
_{I\setminus \gamma J}\omega \right) }{\left\vert J\right\vert ^{\frac{1}{n}}%
}\right) ^{2}\left\Vert \mathsf{P}_{J}^{\limfunc{subgood},\sigma }\mathbf{x}%
\right\Vert _{L^{2}\left( \sigma \right) }^{2}=\sum_{J\in \mathcal{M}_{%
\limfunc{deep}}}\left( \frac{\mathrm{P}^{\alpha }\left( J,\mathbf{1}%
_{I\setminus \gamma J}\omega \right) }{\left\vert J\right\vert ^{\frac{1}{n}}%
}\right) ^{2}\left\Vert \mathsf{P}_{J}^{\limfunc{subgood},\sigma }\mathbf{x}%
\right\Vert _{L^{2}\left( \sigma \right) }^{2}\ .
\end{equation*}%
Let $2<\gamma ^{\prime }<\gamma $ where both $\gamma ^{\prime }$ and $\frac{%
\gamma }{\gamma ^{\prime }}$ will be taken sufficiently large for the
arguments below to be valid - see both (\ref{provided}) and (\ref{provided'}%
) below. For example taking $\gamma ^{\prime }=\sqrt{\gamma }$ and $\gamma
\gg \left( n-\alpha \right) ^{-2}$ works, but is far from optimal. We will
consider the cases $\gamma ^{\prime }J\cap L=\emptyset $ and $\gamma
^{\prime }J\cap L\neq \emptyset $ separately.

Suppose $\gamma ^{\prime }J\cap L=\emptyset $. There is $c>0$ and a finite
sequence $\left\{ \xi _{k}\right\} _{k=1}^{N}$ in $\mathbb{S}^{n-1}$
(actually of the form $\xi _{k}=\left( 0,\xi _{k}^{2},...,\xi
_{k}^{n}\right) $) with the following property. For each $J\in \mathcal{M}_{%
\limfunc{deep}}$ with $\gamma ^{\prime }J\cap L=\emptyset $, there is $1\leq
k=k\left( J\right) \leq N$ such that for $y\in J$ and $x\in I\cap L$, the
linear combination $\xi _{k}\cdot \mathbf{K}^{\alpha ,n}\left( y,x\right) $
is positive and satisfies 
\begin{equation*}
\xi _{k}\cdot \mathbf{K}^{\alpha ,n}\left( y,x\right) =\frac{\xi _{k}\cdot
\left( y-x\right) }{\left\vert y-x\right\vert ^{n+1-\alpha }}\geq c\frac{%
\left\vert J\right\vert ^{\frac{1}{n}}}{\left\vert y-x\right\vert
^{n+1-\alpha }}.
\end{equation*}%
For example, in the plane $n=2$, if $J$ lies above the $x_{1}$-axis $L$,
then for $y\in J$ and $x\in L$ we have $y_{2}\geq \left( \gamma ^{\prime
}-1\right) \left\vert J\right\vert ^{\frac{1}{n}}>\left\vert J\right\vert ^{%
\frac{1}{n}}$ and $x_{2}=0$, hence the estimate 
\begin{equation*}
\left( 0,1\right) \cdot \mathbf{K}^{\alpha ,n}\left( y,x\right) =\frac{%
y_{2}-x_{2}}{\left\vert y-x\right\vert ^{n+1-\alpha }}\geq \frac{\left\vert
J\right\vert ^{\frac{1}{n}}}{\left\vert y-x\right\vert ^{n+1-\alpha }}.
\end{equation*}%
For $J$ below\thinspace $L$ we take the unit vector $\left( 0,-1\right) $ in
place of $\left( 0,1\right) $. Thus for $y\in J\in \mathcal{M}_{\limfunc{deep%
}}$ and $k=k\left( J\right) $ we have the following `weak reversal' of
energy, 
\begin{eqnarray}
\left\vert \mathbf{R}^{\alpha ,n}\left( \mathbf{1}_{I\cap L}\omega \right)
\left( y\right) \right\vert &=&\left\vert \int_{I\cap L}\mathbf{K}^{\alpha
,n}\left( y,x\right) d\omega \left( x\right) \right\vert
\label{weak control} \\
&\geq &\left\vert \int_{I\cap L}\xi _{k}\cdot \mathbf{K}^{\alpha ,n}\left(
y,x\right) d\omega \left( x\right) \right\vert  \notag \\
&\geq &c\int_{I\cap L}\frac{\left\vert J\right\vert ^{\frac{1}{n}}}{%
\left\vert y-x\right\vert ^{n+1-\alpha }}d\omega \left( x\right) \approx c%
\mathrm{P}^{\alpha }\left( J,\mathbf{1}_{I}\omega \right) .  \notag
\end{eqnarray}%
Thus from (\ref{weak control}) and the pairwise disjointedness of $J\in 
\mathcal{M}_{\limfunc{deep}}$, we have%
\begin{eqnarray*}
&&\sum_{\substack{ J\in \mathcal{M}_{\limfunc{deep}}  \\ \gamma ^{\prime
}J\cap L=\emptyset }}\left( \frac{\mathrm{P}^{\alpha }\left( J,\mathbf{1}%
_{I}\omega \right) }{\left\vert J\right\vert ^{\frac{1}{n}}}\right)
^{2}\left\Vert \mathsf{P}_{J}^{\limfunc{subgood},\sigma }\mathbf{x}%
\right\Vert _{L^{2}\left( \sigma \right) }^{2}\leq \sum_{\substack{ J\in 
\mathcal{M}_{\limfunc{deep}}  \\ \gamma ^{\prime }J\cap L=\emptyset }}%
\mathrm{P}^{\alpha }\left( J,\mathbf{1}_{I}\omega \right) ^{2}\left\vert
J\right\vert _{\sigma } \\
&\lesssim &\sum_{J\in \mathcal{M}_{\limfunc{deep}}}\int_{J}\left\vert 
\mathbf{R}^{\alpha ,n}\left( \mathbf{1}_{I\cap L}\omega \right) \left(
y\right) \right\vert ^{2}d\sigma \left( y\right) \\
&\leq &\int_{I}\left\vert \mathbf{R}^{\alpha ,n}\left( \mathbf{1}_{I}\omega
\right) \left( y\right) \right\vert ^{2}d\sigma \left( y\right) \leq \left( 
\mathfrak{T}_{\mathbf{R}^{\alpha ,n}}^{\ast }\right) ^{2}\left\vert
I\right\vert _{\omega }\ .
\end{eqnarray*}

Now we turn to estimating the sum over those cubes $J\in \mathcal{M}_{%
\limfunc{deep}}$ for which $\gamma ^{\prime }J\cap L\neq \emptyset $. In
this case we use the one-dimensional nature of $\omega $ to obtain a strong
reversal of one of the partial energies. Recall the Hilbert transform
inequality for intervals $J$ and $I$ with $2J\subset I$ and $\limfunc{supp}%
\mu \subset \mathbb{R}\setminus I$: 
\begin{eqnarray}
\sup_{y,z\in J}\frac{H\mu \left( y\right) -H\mu \left( z\right) }{y-z}
&=&\int_{\mathbb{R}\setminus I}\left\{ \frac{\frac{1}{x-y}-\frac{1}{x-z}}{y-z%
}\right\} d\mu \left( x\right)  \label{Hilbert} \\
&=&\int_{\mathbb{R}\setminus I}\frac{1}{\left( x-y\right) \left( x-z\right) }%
d\mu \left( x\right) \approx \frac{\mathrm{P}\left( J,\mu \right) }{%
\left\vert J\right\vert }.  \notag
\end{eqnarray}%
We wish to obtain a similar control in the situation at hand, but the matter
is now complicated by the extra dimensions. Fix $y=\left( y^{1},y^{\prime
}\right) ,z=\left( z^{1},z^{\prime }\right) \in J$ and $x=\left(
x^{1},0\right) \in L\setminus \gamma J$. We consider first the case%
\begin{equation}
\left\vert y^{\prime }-z^{\prime }\right\vert \leq C_{0}\left\vert
y^{1}-z^{1}\right\vert ,  \label{case 1}
\end{equation}%
where $C_{0}$ is a positive constant satisfying (\ref{provided''}) below.
Now the first component $R_{1}^{\alpha ,n}$ is `positive' in the direction
of the $x^{1}$-axis $L$, and so for $\left( y^{1},y^{\prime }\right) ,\left(
z^{1},z^{\prime }\right) \in J$, we write%
\begin{eqnarray*}
&&\frac{R_{1}^{\alpha ,n}\mathbf{1}_{I\setminus \gamma J}\omega \left(
y^{1},y^{\prime }\right) -R_{1}^{\alpha ,n}\mathbf{1}_{I\setminus \gamma
J}\omega \left( z^{1},z^{\prime }\right) }{y^{1}-z^{1}} \\
&=&\int_{I\setminus \gamma J}\left\{ \frac{K_{1}^{\alpha ,n}\left( \left(
y^{1},y^{\prime }\right) ,x\right) -K_{1}^{\alpha ,n}\left( \left(
z^{1},z^{\prime }\right) ,x\right) }{y^{1}-z^{1}}\right\} d\omega \left(
x\right) \\
&=&\int_{I\setminus \gamma J}\left\{ \frac{\frac{y^{1}-x^{1}}{\left\vert
y-x\right\vert ^{n+1-\alpha }}-\frac{z^{1}-x^{1}}{\left\vert z-x\right\vert
^{n+1-\alpha }}}{y^{1}-z^{1}}\right\} d\omega \left( x\right) \ .
\end{eqnarray*}%
For $0\leq t\leq 1$ define 
\begin{eqnarray*}
w_{t} &\equiv &ty+\left( 1-t\right) z=z+t\left( y-z\right) , \\
w_{t}-x &=&t\left( y-x\right) +\left( 1-t\right) \left( z-x\right) ,
\end{eqnarray*}%
and 
\begin{equation*}
\Phi \left( t\right) \equiv \frac{w_{t}^{1}-x^{1}}{\left\vert
w_{t}-x\right\vert ^{n+1-\alpha }}\ ,
\end{equation*}%
so that%
\begin{equation*}
\frac{y^{1}-x^{1}}{\left\vert y-x\right\vert ^{n+1-\alpha }}-\frac{%
z^{1}-x^{1}}{\left\vert z-x\right\vert ^{n+1-\alpha }}=\Phi \left( 1\right)
-\Phi \left( 0\right) =\int_{0}^{1}\Phi ^{\prime }\left( t\right) dt\ .
\end{equation*}%
Then using $\nabla \left\vert \xi \right\vert ^{\tau }=\tau \left\vert \xi
\right\vert ^{\tau -2}\xi $ we compute that%
\begin{eqnarray*}
\frac{d}{dt}\Phi \left( t\right) &=&\frac{\left( y^{1}-z^{1}\right) }{%
\left\vert w_{t}-x\right\vert ^{n+1-\alpha }}+\left( w_{t}^{1}-x^{1}\right)
\left( y-z\right) \cdot \nabla \left\vert w_{t}-x\right\vert ^{-n-1+\alpha }
\\
&=&\frac{\left( y^{1}-z^{1}\right) }{\left\vert w_{t}-x\right\vert
^{n+1-\alpha }}-\left( n+1-\alpha \right) \left( w_{t}^{1}-x^{1}\right) 
\frac{\left( w_{t}-x\right) \cdot \left( y-z\right) }{\left\vert
w_{t}-x\right\vert ^{n+3-\alpha }} \\
&=&\frac{\left( y^{1}-z^{1}\right) }{\left\vert w_{t}-x\right\vert
^{n+1-\alpha }}-\left( n+1-\alpha \right) \left( w_{t}^{1}-x^{1}\right) 
\frac{\left( w_{t}^{1}-x^{1}\right) \left( y^{1}-z^{1}\right) }{\left\vert
w_{t}-x\right\vert ^{n+3-\alpha }} \\
&&-\left( n+1-\alpha \right) \left( w_{t}^{1}-x^{1}\right) \frac{\left(
w_{t}^{\prime }-x^{\prime }\right) \cdot \left( y^{\prime }-z^{\prime
}\right) }{\left\vert w_{t}-x\right\vert ^{n+3-\alpha }} \\
&=&\left( y^{1}-z^{1}\right) \left\{ \frac{\left\vert w_{t}-x\right\vert ^{2}%
}{\left\vert w_{t}-x\right\vert ^{n+3-\alpha }}-\left( n+1-\alpha \right) 
\frac{\left\vert w_{t}^{1}-x^{1}\right\vert ^{2}}{\left\vert
w_{t}-x\right\vert ^{n+3-\alpha }}\right\} \\
&&+\left( y^{1}-z^{1}\right) \left\{ -\left( n+1-\alpha \right) \left( \frac{%
w_{t}^{1}-x^{1}}{y^{1}-z^{1}}\right) \frac{\left( w_{t}^{\prime }-x^{\prime
}\right) \cdot \left( y^{\prime }-z^{\prime }\right) }{\left\vert
w_{t}-x\right\vert ^{n+3-\alpha }}\right\} \\
&\equiv &\left( y^{1}-z^{1}\right) \left\{ A\left( t\right) +B\left(
t\right) \right\} .
\end{eqnarray*}%
Now $\left\vert w_{t}^{1}-x^{1}\right\vert \approx \left\vert y-x\right\vert 
$ and $\left\vert w_{t}^{\prime }-x^{\prime }\right\vert =\left\vert
w_{t}^{\prime }\right\vert \approx \left\vert y^{\prime }\right\vert \leq
\gamma ^{\prime }\frac{\left\vert y-x\right\vert }{\gamma }$ because $\gamma
^{\prime }J\cap L=\emptyset $, and so if $\gamma \gg \gamma ^{\prime }$ we
obtain using$\left\vert y-x\right\vert \approx \left\vert
w_{t}^{1}-x^{1}\right\vert $ that%
\begin{equation*}
\left\vert w_{t}^{\prime }-x^{\prime }\right\vert \leq \sqrt{\frac{n-\alpha 
}{2}}\left\vert w_{t}^{1}-x^{1}\right\vert ,
\end{equation*}%
and hence that%
\begin{eqnarray*}
-A\left( t\right) &=&-\frac{\left\vert w_{t}-x\right\vert ^{2}}{\left\vert
w_{t}-x\right\vert ^{n+3-\alpha }}+\left( n+1-\alpha \right) \frac{%
\left\vert w_{t}^{1}-x^{1}\right\vert ^{2}}{\left\vert w_{t}-x\right\vert
^{n+3-\alpha }} \\
&=&\frac{-\left\vert w_{t}-x\right\vert ^{2}+\left( n+1-\alpha \right)
\left\vert w_{t}^{1}-x^{1}\right\vert ^{2}}{\left\vert w_{t}-x\right\vert
^{n+3-\alpha }} \\
&=&\frac{-\left\vert w_{t}^{\prime }-x^{\prime }\right\vert ^{2}+\left(
n-\alpha \right) \left( w_{t}^{1}-x^{1}\right) ^{2}}{\left\vert
w_{t}-x\right\vert ^{n+3-\alpha }} \\
&\approx &\left( n-\alpha \right) \frac{\left( w_{t}^{1}-x^{1}\right) ^{2}}{%
\left\vert w_{t}-x\right\vert ^{n+3-\alpha }}\ .
\end{eqnarray*}%
Now from our assumption (\ref{case 1}) we have%
\begin{eqnarray*}
\left\vert B\left( t\right) \right\vert &=&\left\vert \left( n+1-\alpha
\right) \left( \frac{w_{t}^{1}-x^{1}}{y^{1}-z^{1}}\right) \frac{\left(
w_{t}^{\prime }-x^{\prime }\right) \cdot \left( y^{\prime }-z^{\prime
}\right) }{\left\vert w_{t}-x\right\vert ^{n+3-\alpha }}\right\vert \\
&\leq &\left( n+1-\alpha \right) \frac{\left\vert w_{t}^{1}-x^{1}\right\vert 
}{\left\vert y^{1}-z^{1}\right\vert }\frac{\left\vert w_{t}^{\prime
}-x^{\prime }\right\vert \left\vert y^{\prime }-z^{\prime }\right\vert }{%
\left\vert w_{t}-x\right\vert ^{n+3-\alpha }} \\
&\leq &C_{0}\left( n+1-\alpha \right) \frac{\left\vert
w_{t}^{1}-x^{1}\right\vert \left\vert w_{t}^{\prime }-x^{\prime }\right\vert 
}{\left\vert w_{t}-x\right\vert ^{n+3-\alpha }} \\
&\leq &C_{0}\left( n+1-\alpha \right) \frac{\gamma ^{\prime }}{\gamma }\frac{%
\left\vert w_{t}^{1}-x^{1}\right\vert ^{2}}{\left\vert w_{t}-x\right\vert
^{n+3-\alpha }}\ll \frac{1}{2}\left( n-\alpha \right) \frac{\left(
w_{t}^{1}-x^{1}\right) ^{2}}{\left\vert w_{t}-x\right\vert ^{n+3-\alpha }}
\end{eqnarray*}%
if 
\begin{equation}
C_{0}\ll \frac{\gamma }{2\gamma ^{\prime }}\frac{n-\alpha }{n+1-\alpha }.
\label{provided''}
\end{equation}%
Thus altogether in case (\ref{case 1}) we have 
\begin{eqnarray*}
&&\left\vert R_{1}^{\alpha ,n}\mathbf{1}_{I\setminus \gamma J}\omega \left(
y^{1},y^{\prime }\right) -R_{1}^{\alpha ,n}\mathbf{1}_{I\setminus \gamma
J}\omega \left( z^{1},z^{\prime }\right) \right\vert \\
&\approx &\left\vert y^{1}-z^{1}\right\vert \left\vert \int_{I\setminus
\gamma J}\frac{\int_{0}^{1}\frac{d}{dt}\Phi \left( t\right) dt}{y^{1}-z^{1}}%
d\omega \left( x\right) \right\vert \\
&\approx &\left\vert y^{1}-z^{1}\right\vert \left\vert \int_{I\setminus
\gamma J}\int_{0}^{1}\left\{ A\left( t\right) +B\left( t\right) \right\}
dtd\omega \left( x\right) \right\vert \\
&\approx &\left\vert y^{1}-z^{1}\right\vert \left\vert \int_{I\setminus
\gamma J}\int_{0}^{1}\left\{ \left( n-\alpha \right) \frac{\left(
w_{t}^{1}-x^{1}\right) ^{2}}{\left\vert w_{t}-x\right\vert ^{n+3-\alpha }}%
\right\} dtd\omega \left( x\right) \right\vert \\
&\approx &\left\vert y^{1}-z^{1}\right\vert \int_{I\setminus \gamma J}\frac{%
\left( c_{J}^{1}-x^{1}\right) ^{2}}{\left\vert c_{J}-x\right\vert
^{n+3-\alpha }}d\omega \left( x\right) \\
&\approx &\left\vert y^{1}-z^{1}\right\vert \frac{\mathrm{P}^{\alpha }\left(
J,\mathbf{1}_{I\setminus \gamma J}\omega \right) }{\left\vert J\right\vert ^{%
\frac{1}{n}}}\ .
\end{eqnarray*}

On the other hand, in the case that%
\begin{equation}
\left\vert y^{\prime }-z^{\prime }\right\vert >C_{0}\left\vert
y^{1}-z^{1}\right\vert ,  \label{case 2}
\end{equation}%
we write%
\begin{eqnarray*}
\left( \mathbf{R}^{\alpha ,n}\right) ^{\prime } &=&\left( R_{2}^{\alpha
,n},...,R_{n}^{\alpha ,n}\right) , \\
\mathbf{\Phi }\left( t\right) &=&\frac{w_{t}^{\prime }-x^{\prime }}{%
\left\vert w_{t}-x\right\vert ^{n+1-\alpha }},
\end{eqnarray*}%
with $w_{t}=ty+\left( 1-t\right) z$ as before. Then as above we obtain 
\begin{equation*}
\frac{y^{\prime }-x^{\prime }}{\left\vert y-x\right\vert ^{n+1-\alpha }}-%
\frac{z^{\prime }-x^{\prime }}{\left\vert z-x\right\vert ^{n+1-\alpha }}=%
\mathbf{\Phi }\left( 1\right) -\mathbf{\Phi }\left( 0\right) =\int_{0}^{1}%
\frac{d}{dt}\mathbf{\Phi }\left( t\right) dt,
\end{equation*}%
where if we write $\widehat{y^{k}}\equiv \left(
y^{1},...,y^{k-1},0,y^{k+1},...,y^{n}\right) $, we have%
\begin{eqnarray*}
\frac{d}{dt}\mathbf{\Phi }\left( t\right) &=&\left\{ \frac{d}{dt}\Phi
_{k}\left( t\right) \right\} _{k=2}^{n} \\
&=&\left\{ \left( y^{k}-z^{k}\right) \left[ \frac{\left\vert
w_{t}-x\right\vert ^{2}}{\left\vert w_{t}-x\right\vert ^{n+3-\alpha }}%
-\left( n+1-\alpha \right) \frac{\left\vert w_{t}^{k}-x^{k}\right\vert ^{2}}{%
\left\vert w_{t}-x\right\vert ^{n+3-\alpha }}\right] \right\} _{k=2}^{n} \\
&&-\left\{ \left( n+1-\alpha \right) \left( w_{t}^{k}-x^{k}\right) \frac{%
\left( \widehat{w_{t}^{k}}-\widehat{x^{k}}\right) \cdot \left( \widehat{y^{k}%
}-\widehat{z^{k}}\right) }{\left\vert w_{t}-x\right\vert ^{n+3-\alpha }}%
\right\} _{k=2}^{n} \\
&\equiv &\left\{ \left( y^{k}-z^{k}\right) A_{k}\left( t\right) \right\}
_{k=2}^{n}+\left\{ V_{k}\left( t\right) \right\} _{k=2}^{n}\equiv \mathbf{U}%
\left( t\right) +\mathbf{V}\left( t\right) \ .
\end{eqnarray*}%
Now for $2\leq k\leq n$ we have $x^{k}=0$ and so%
\begin{eqnarray*}
A_{k}\left( t\right) &=&\frac{\left\vert w_{t}-x\right\vert ^{2}}{\left\vert
w_{t}-x\right\vert ^{n+3-\alpha }}-\left( n+1-\alpha \right) \frac{%
\left\vert w_{t}^{k}\right\vert ^{2}}{\left\vert w_{t}-x\right\vert
^{n+3-\alpha }} \\
&=&\frac{\left\vert w_{t}-x\right\vert ^{2}-\left( n+1-\alpha \right)
\left\vert w_{t}^{k}\right\vert ^{2}}{\left\vert w_{t}-x\right\vert
^{n+3-\alpha }} \\
&=&\frac{\left\vert w_{t}^{1}-x^{1}\right\vert ^{2}+\sum_{j\neq
1,k}\left\vert w_{t}^{j}\right\vert ^{2}-\left( n-\alpha \right) \left(
w_{t}^{k}\right) ^{2}}{\left\vert w_{t}-x\right\vert ^{n+3-\alpha }} \\
&\approx &\left( n-\alpha \right) \frac{\left\vert
w_{t}^{1}-x^{1}\right\vert ^{2}}{\left\vert w_{t}-x\right\vert ^{n+3-\alpha }%
}\approx \left( n-\alpha \right) \frac{1}{\left\vert c_{J}-x\right\vert
^{n+1-\alpha }}\ .
\end{eqnarray*}%
Thus we have%
\begin{equation*}
\int_{I\setminus \gamma J}A_{k}\left( t\right) d\omega \left( x\right)
\approx \left( n-\alpha \right) \int_{I\setminus \gamma J}\frac{1}{%
\left\vert c_{J}-x\right\vert ^{n+1-\alpha }}d\omega \left( x\right) \approx
\left( n-\alpha \right) \frac{\mathrm{P}^{\alpha }\left( J,\mathbf{1}%
_{I\setminus \gamma J}\omega \right) }{\left\vert J\right\vert ^{\frac{1}{n}}%
},
\end{equation*}%
and hence%
\begin{eqnarray*}
\left\vert \int_{I\setminus \gamma J}\int_{0}^{1}\mathbf{U}\left( t\right)
dtd\omega \left( x\right) \right\vert ^{2} &=&\left\vert \int_{I\setminus
\gamma J}\int_{0}^{1}\left\{ \left( y^{k}-z^{k}\right) A_{k}\left( t\right)
\right\} _{k=2}^{n}dtd\omega \left( x\right) \right\vert ^{2} \\
&=&\sum_{k=2}^{n}\left( y^{k}-z^{k}\right) ^{2}\left\vert \int_{I\setminus
\gamma J}\int_{0}^{1}\left\{ A_{k}\left( t\right) \right\}
_{k=2}^{n}dtd\omega \left( x\right) \right\vert ^{2} \\
&\approx &\sum_{k=2}^{n}\left( y^{k}-z^{k}\right) ^{2}\left( n-\alpha
\right) ^{2}\left( \frac{\mathrm{P}^{\alpha }\left( J,\mathbf{1}_{I\setminus
\gamma J}\omega \right) }{\left\vert J\right\vert ^{\frac{1}{n}}}\right) ^{2}
\\
&\approx &\left( n-\alpha \right) ^{2}\left\vert y^{\prime }-z^{\prime
}\right\vert ^{2}\left( \frac{\mathrm{P}^{\alpha }\left( J,\mathbf{1}%
_{I\setminus \gamma J}\omega \right) }{\left\vert J\right\vert ^{\frac{1}{n}}%
}\right) ^{2}.
\end{eqnarray*}%
For $2\leq k\leq n$ we also have using (\ref{case 2}) that%
\begin{eqnarray*}
\left\vert V_{k}\left( t\right) \right\vert &=&\left\vert \left( n+1-\alpha
\right) \left( w_{t}^{k}-x^{k}\right) \frac{\left( \widehat{w_{t}^{k}}-%
\widehat{x^{k}}\right) \cdot \left( \widehat{y^{k}}-\widehat{z^{k}}\right) }{%
\left\vert w_{t}-x\right\vert ^{n+3-\alpha }}\right\vert \\
&\leq &\left( n+1-\alpha \right) \left\vert w_{t}^{k}\right\vert \frac{%
\left\vert w_{t}^{1}-x^{1}\right\vert \left\vert y^{1}-z^{1}\right\vert
+\sum_{j\neq 1,k}\left\vert w_{t}^{j}\right\vert \left\vert
y^{j}-z^{j}\right\vert }{\left\vert w_{t}-x\right\vert ^{n+3-\alpha }} \\
&\leq &\left( n+1-\alpha \right) \left\{ \frac{\left\vert
w_{t}^{k}\right\vert \left\vert y^{1}-z^{1}\right\vert }{\left\vert
w_{t}-x\right\vert ^{n+2-\alpha }}+\sum_{j\neq 1,k}\frac{\left\vert
w_{t}^{k}\right\vert \left\vert w_{t}^{j}\right\vert \left\vert
y^{j}-z^{j}\right\vert }{\left\vert w_{t}-x\right\vert ^{n+3-\alpha }}%
\right\} \\
&\leq &\left( n+1-\alpha \right) \left\{ \frac{\gamma ^{\prime }}{\gamma }%
\frac{\left\vert y^{1}-z^{1}\right\vert }{\left\vert w_{t}-x\right\vert
^{n+1-\alpha }}+\left( \gamma ^{\prime }\right) ^{2}\frac{\left\vert
y^{\prime }-z^{\prime }\right\vert }{\left\vert w_{t}-x\right\vert
^{n+3-\alpha }}\right\} \\
&\lesssim &\left( n+1-\alpha \right) \left\{ \frac{\gamma ^{\prime }}{\gamma
C_{0}}\frac{\left\vert y^{\prime }-z^{\prime }\right\vert }{\left\vert
c_{J}-x\right\vert ^{n+1-\alpha }}+\left( \frac{\gamma ^{\prime }}{\gamma }%
\right) ^{2}\frac{\left\vert y^{\prime }-z^{\prime }\right\vert }{\left\vert
c_{J}-x\right\vert ^{n+1-\alpha }}\right\} .
\end{eqnarray*}%
Thus%
\begin{eqnarray*}
\left\vert \int_{I\setminus \gamma J}\int_{0}^{1}\mathbf{V}\left( t\right)
dtd\omega \left( x\right) \right\vert &\lesssim &\left( n+1-\alpha \right)
\left\{ \frac{\gamma ^{\prime }}{\gamma C_{0}}+\left( \frac{\gamma ^{\prime }%
}{\gamma }\right) ^{2}\right\} \int_{I\setminus \gamma J}\frac{\left\vert
y^{\prime }-z^{\prime }\right\vert }{\left\vert c_{J}-x\right\vert
^{n+1-\alpha }}d\omega \left( x\right) \\
&\lesssim &\left( n+1-\alpha \right) \left\{ \frac{\gamma ^{\prime }}{\gamma
C_{0}}+\left( \frac{\gamma ^{\prime }}{\gamma }\right) ^{2}\right\}
\left\vert y^{\prime }-z^{\prime }\right\vert \frac{\mathrm{P}^{\alpha
}\left( J,\mathbf{1}_{I\setminus \gamma J}\omega \right) }{\left\vert
J\right\vert ^{\frac{1}{n}}},
\end{eqnarray*}%
and so%
\begin{equation*}
\left\vert \int_{I\setminus \gamma J}\int_{0}^{1}\mathbf{V}\left( t\right)
dtd\omega \left( x\right) \right\vert \leq \frac{1}{2}\left\vert
\int_{I\setminus \gamma J}\int_{0}^{1}\mathbf{U}\left( t\right) dtd\omega
\left( x\right) \right\vert \ ,
\end{equation*}%
provided%
\begin{eqnarray}
&&\left( n+1-\alpha \right) \left\{ \frac{\gamma ^{\prime }}{\gamma C_{0}}%
+\left( \frac{\gamma ^{\prime }}{\gamma }\right) ^{2}\right\} \ll n-\alpha ,
\label{provided} \\
&&\text{i.e. }\left( \frac{n+1-\alpha }{n-\alpha }\right) ^{2}\left( \frac{%
\gamma ^{\prime }}{\gamma }\right) ^{2}\ll 1\ ,  \notag
\end{eqnarray}%
where we have used (\ref{provided''}) with an optimal $C_{0}$. Then if both (%
\ref{case 2}) and (\ref{provided}) hold we have 
\begin{eqnarray*}
&&\left\vert \left( \mathbf{R}^{\alpha ,n}\right) ^{\prime }\mathbf{1}%
_{I\setminus \gamma J}\omega \left( y^{1},y^{\prime }\right) -\left( \mathbf{%
R}^{\alpha ,n}\right) ^{\prime }\mathbf{1}_{I\setminus \gamma J}\omega
\left( z^{1},z^{\prime }\right) \right\vert \\
&=&\left\vert \int_{I\setminus \gamma J}\left\{ \frac{y^{k}-x^{k}}{%
\left\vert y-x\right\vert ^{n+1-\alpha }}-\frac{z^{k}-x^{k}}{\left\vert
z-x\right\vert ^{n+1-\alpha }}\right\} _{k=2}^{n}d\omega \left( x\right)
\right\vert \\
&=&\left\vert \int_{I\setminus \gamma J}\int_{0}^{1}\mathbf{\Phi }^{\prime
}\left( t\right) dtd\omega \left( x\right) \right\vert \\
&\geq &\left\vert \int_{I\setminus \gamma J}\int_{0}^{1}\mathbf{U}\left(
t\right) dtd\omega \left( x\right) \right\vert -\left\vert \int_{I\setminus
\gamma J}\int_{0}^{1}\mathbf{V}\left( t\right) dtd\omega \left( x\right)
\right\vert \\
&\geq &\frac{1}{2}\left\vert \int_{I\setminus \gamma J}\int_{0}^{1}\mathbf{U}%
\left( t\right) dtd\omega \left( x\right) \right\vert \\
&\gtrsim &C_{0}\int_{I\setminus \gamma J}\int_{0}^{1}\frac{\left\vert
y^{1}-z^{1}\right\vert }{\left\vert c_{J}-x\right\vert ^{n+1-\alpha }}%
dtd\omega \left( x\right) \approx C_{0}\left\vert y^{1}-z^{1}\right\vert 
\frac{\mathrm{P}^{\alpha }\left( J,\mathbf{1}_{I\setminus \gamma J}\omega
\right) }{\left\vert J\right\vert ^{\frac{1}{n}}}\ .
\end{eqnarray*}

Combining the inequalities from each case (\ref{case 1}) and (\ref{case 2})
above, and assuming (\ref{provided}), we conclude that for all $y,z\in J$ we
have the following `strong reversal' of the $1$-partial\ energy,%
\begin{equation*}
\left\vert y^{1}-z^{1}\right\vert ^{2}\left( \frac{\mathrm{P}^{\alpha
}\left( J,\mathbf{1}_{I\setminus \gamma J}\omega \right) }{\left\vert
J\right\vert ^{\frac{1}{n}}}\right) ^{2}\lesssim \left\vert \mathbf{R}%
^{\alpha ,n}\mathbf{1}_{I\setminus \gamma J}\omega \left( y^{1},y^{\prime
}\right) -\mathbf{R}^{\alpha ,n}\mathbf{1}_{I\setminus \gamma J}\omega
\left( z^{1},z^{\prime }\right) \right\vert ^{2}.
\end{equation*}%
Thus we have 
\begin{eqnarray*}
&&\sum_{\substack{ J\in \mathcal{M}_{\limfunc{deep}}  \\ \gamma ^{\prime
}J\cap L\neq \emptyset }}\left( \frac{\mathrm{P}^{\alpha }\left( J,\mathbf{1}%
_{I\setminus \gamma J}\omega \right) }{\left\vert J\right\vert ^{\frac{1}{n}}%
}\right) ^{2}\int_{J}\left\vert y^{1}-\mathbb{E}_{J}^{\sigma
}y^{1}\right\vert ^{2}d\sigma \left( y\right) \\
&=&\frac{1}{2}\sum_{\substack{ J\in \mathcal{M}_{\limfunc{deep}}  \\ \gamma
^{\prime }J\cap L\neq \emptyset }}\left( \frac{\mathrm{P}^{\alpha }\left( J,%
\mathbf{1}_{I\setminus \gamma J}\omega \right) }{\left\vert J\right\vert ^{%
\frac{1}{n}}}\right) ^{2}\frac{1}{\left\vert J\right\vert _{\sigma }}%
\int_{J}\int_{J}\left( y^{1}-z^{1}\right) ^{2}d\sigma \left( y\right)
d\sigma \left( z\right) \\
&\lesssim &\sum_{\substack{ J\in \mathcal{M}_{\limfunc{deep}}  \\ \gamma
^{\prime }J\cap L\neq \emptyset }}\frac{1}{\left\vert J\right\vert _{\sigma }%
}\int_{J}\int_{J}\left\vert \mathbf{R}^{\alpha ,n}\mathbf{1}_{I\setminus
\gamma J}\omega \left( y^{1},y^{\prime }\right) -\mathbf{R}^{\alpha ,n}%
\mathbf{1}_{I\setminus \gamma J}\omega \left( z^{1},z^{\prime }\right)
\right\vert ^{2}d\sigma \left( y\right) d\sigma \left( z\right) \\
&\lesssim &\sum_{\substack{ J\in \mathcal{M}_{\limfunc{deep}}  \\ \gamma
^{\prime }J\cap L\neq \emptyset }}\int_{J}\left\vert \mathbf{R}^{\alpha ,n}%
\mathbf{1}_{I\setminus \gamma J}\omega \left( y^{1},y^{\prime }\right)
\right\vert ^{2}d\sigma \left( y\right) \\
&\leq &\sum_{J\in \mathcal{M}_{\limfunc{deep}}}\int_{J}\left\vert \mathbf{R}%
^{\alpha ,n}\mathbf{1}_{I}\omega \left( y^{1},y^{\prime }\right) \right\vert
^{2}d\sigma \left( y\right) +\sum_{J\in \mathcal{M}_{\limfunc{deep}%
}}\int_{J}\left\vert \mathbf{R}^{\alpha ,n}\mathbf{1}_{\gamma J}\omega
\left( y^{1},y^{\prime }\right) \right\vert ^{2}d\sigma \left( y\right) ,
\end{eqnarray*}%
and now we obtain in the usual way that this is bounded by%
\begin{eqnarray*}
&&\int_{I}\left\vert \mathbf{R}^{\alpha ,n}\mathbf{1}_{I}\omega \left(
y^{1},y^{\prime }\right) \right\vert ^{2}d\sigma \left( y\right) +\sum_{J\in 
\mathcal{M}}\left( \mathfrak{T}_{\mathbf{R}^{\alpha ,n}}^{\ast }\right)
^{2}\left\vert \gamma J\right\vert _{\omega } \\
&\leq &\left( \mathfrak{T}_{\mathbf{R}^{\alpha ,n}}^{\ast }\right)
^{2}\left\vert I\right\vert _{\omega }+\beta \left( \mathfrak{T}_{\mathbf{R}%
^{\alpha ,n}}^{\ast }\right) ^{2}\left\vert I\right\vert _{\omega }\lesssim
\left( \mathfrak{T}_{\mathbf{R}^{\alpha ,n}}^{\ast }\right) ^{2}\left\vert
I\right\vert _{\omega }\ .
\end{eqnarray*}

Now we turn to the other partial energies and begin with the estimate that
for $2\leq j\leq n$, we have the following `weak reversal' of energy,%
\begin{eqnarray*}
\left\vert R_{j}^{\alpha ,n}\mathbf{1}_{I\setminus \gamma J}\omega \left(
y\right) \right\vert &=&\left\vert \int_{I\setminus \gamma J}\frac{y^{j}-0}{%
\left\vert y-x\right\vert ^{n+1-\alpha }}d\omega \left( x_{1},0...,0\right)
\right\vert \\
&\approx &\left\vert \frac{y^{j}}{\left\vert J\right\vert ^{\frac{1}{n}}}%
\int_{I\setminus \gamma J}\frac{\left\vert J\right\vert ^{\frac{1}{n}}}{%
\left\vert y-x\right\vert ^{n+1-\alpha }}d\omega \left( x_{1},0...,0\right)
\right\vert \\
&\approx &\left\vert y^{j}\right\vert \frac{\mathrm{P}^{\alpha }\left( J,%
\mathbf{1}_{I\setminus \gamma J}\omega \right) }{\left\vert J\right\vert ^{%
\frac{1}{n}}}.
\end{eqnarray*}%
Thus for $2\leq j\leq n$, we use $\int_{J}\left\vert y^{j}-\mathbb{E}%
_{J}^{\sigma }y^{j}\right\vert ^{2}d\sigma \left( y\right) \leq
\int_{J}\left\vert y^{j}\right\vert ^{2}d\sigma \left( y\right) $ to obtain%
\begin{eqnarray*}
&&\sum_{\substack{ J\in \mathcal{M}_{\limfunc{deep}}  \\ \gamma ^{\prime
}J\cap L\neq \emptyset }}\left( \frac{\mathrm{P}^{\alpha }\left( J,\mathbf{1}%
_{I\setminus \gamma J}\omega \right) }{\left\vert J\right\vert ^{\frac{1}{n}}%
}\right) ^{2}\int_{J}\left\vert y^{j}-\mathbb{E}_{J}^{\sigma
}y^{j}\right\vert ^{2}d\sigma \left( y\right) \\
&\leq &\sum_{\substack{ J\in \mathcal{M}_{\limfunc{deep}}  \\ \gamma
^{\prime }J\cap L\neq \emptyset }}\left( \frac{\mathrm{P}^{\alpha }\left( J,%
\mathbf{1}_{I\setminus \gamma J}\omega \right) }{\left\vert J\right\vert ^{%
\frac{1}{n}}}\right) ^{2}\int_{J}\left\vert y^{j}\right\vert ^{2}d\sigma
\left( y\right) =\sum_{\substack{ J\in \mathcal{M}_{\limfunc{deep}}  \\ %
\gamma ^{\prime }J\cap L\neq \emptyset }}\int_{J}\left( \frac{\mathrm{P}%
^{\alpha }\left( J,\mathbf{1}_{I\setminus \gamma J}\omega \right) }{%
\left\vert J\right\vert ^{\frac{1}{n}}}\right) ^{2}\left\vert
y^{j}\right\vert ^{2}d\sigma \left( y\right) \\
&\lesssim &\sum_{\substack{ J\in \mathcal{M}_{\limfunc{deep}}  \\ \gamma
^{\prime }J\cap L\neq \emptyset }}\int_{J}\left\vert R_{j}^{\alpha ,n}%
\mathbf{1}_{I\setminus \gamma J}\omega \left( y\right) \right\vert
^{2}d\sigma \left( y\right) \lesssim \int_{I}\left\vert \mathbf{R}^{\alpha
,n}\mathbf{1}_{I}\omega \left( y\right) \right\vert ^{2}d\sigma \left(
y\right) \leq \left( \mathfrak{T}_{\mathbf{R}^{\alpha ,n}}^{\ast }\right)
^{2}\left\vert I\right\vert _{\omega }\ .
\end{eqnarray*}%
Summing these estimates for $j=1$ and $2\leq j\leq n$ completes the proof of
the dual energy condition $\mathcal{E}_{\alpha }^{\ast }\lesssim \mathfrak{T}%
_{\mathbf{R}^{\alpha ,n}}^{\ast }+\sqrt{\mathcal{A}_{2}^{\alpha ,\ast }}$.

\subsection{Forward energy condition}

Now we turn to proving the (forward) energy condition $\mathcal{E}_{\alpha
}\lesssim \mathfrak{T}_{\mathbf{T}^{\alpha ,n}}+\sqrt{\mathcal{A}%
_{2}^{\alpha }}$. We must show%
\begin{equation*}
\sup_{\ell \geq 0}\sum_{r=1}^{\infty }\sum_{J\in \mathcal{M}_{\limfunc{deep}%
}^{\ell }\left( I_{r}\right) }\left( \frac{\mathrm{P}^{\alpha }\left( J,%
\mathbf{1}_{I\setminus J^{\ast }}\sigma \right) }{\left\vert J\right\vert ^{%
\frac{1}{n}}}\right) ^{2}\left\Vert \mathsf{P}_{J}^{\limfunc{subgood},\omega
}\mathbf{x}\right\Vert _{L^{2}\left( \omega \right) }^{2}\leq \left( 
\mathfrak{T}_{\mathbf{R}^{\alpha ,n}}^{2}+\mathcal{A}_{2}^{\alpha }\right)
\left\vert I\right\vert _{\sigma }\ ,
\end{equation*}%
for all partitions of a dyadic cube $I=\overset{\cdot }{\dbigcup\limits_{r%
\geq 1}}I_{r}$ into dyadic subcubes $I_{r}$. We again fix $\ell \geq 0$ and
suppress both $\ell $ and $\mathbf{r}$ in the notation $\mathcal{M}_{%
\limfunc{deep}}\left( I_{r}\right) =\mathcal{M}_{\mathbf{r}-\limfunc{deep}%
}^{\ell }\left( I_{r}\right) $. We may assume that all the cubes $J$
intersect $\limfunc{supp}\omega $, hence that all the cubes $I_{r}$ and $J$
intersect $L$, which contains $\limfunc{supp}\omega $. Let $\mathcal{I}%
_{r}=I_{r}\cap L$ and $\mathcal{J}=J\cap L$ for these cubes. We must show%
\begin{equation*}
\sum_{r=1}^{\infty }\sum_{J\in \mathcal{M}_{\limfunc{deep}}\left(
I_{r}\right) }\left( \frac{\mathrm{P}^{\alpha }\left( J,\mathbf{1}%
_{I\setminus J^{\ast }}\sigma \right) }{\left\vert J\right\vert ^{\frac{1}{n}%
}}\right) ^{2}\left\Vert \mathsf{P}_{J}^{\limfunc{subgood},\omega }\mathbf{x}%
\right\Vert _{L^{2}\left( \omega \right) }^{2}\leq \left( \mathfrak{T}_{%
\mathbf{R}^{\alpha ,n}}^{2}+\mathcal{A}_{2}^{\alpha }\right) \left\vert
I\right\vert _{\sigma }\ .
\end{equation*}%
Let $\mathcal{M}_{\limfunc{deep}}=\dbigcup\limits_{r=1}^{\infty }\mathcal{M}%
_{\limfunc{deep}}\left( I_{r}\right) $ as above, and for each $J\in \mathcal{%
M}_{\limfunc{deep}}$, make the decomposition%
\begin{equation*}
I\setminus J^{\ast }=\text{\textsc{E}}\left( J^{\ast }\right) \dot{\cup}%
\text{\textsc{S}}\left( J^{\ast }\right)
\end{equation*}%
of $I\setminus J^{\ast }$ into \emph{end} $E\left( J^{\ast }\right) $ and 
\emph{side} $S\left( J^{\ast }\right) $ disjoint pieces defined by%
\begin{eqnarray*}
\text{\textsc{E}}\left( J^{\ast }\right) &\equiv &I\cap \left\{ \left(
y^{1},y^{\prime }\right) :\left\vert y^{1}-c_{J}^{1}\right\vert \geq \frac{%
\gamma }{2}\left\vert J\right\vert ^{\frac{1}{n}}\text{ and }\left\vert
y^{\prime }-c_{J}^{\prime }\right\vert \leq \frac{1}{\gamma }\left\vert
y^{1}-c_{J}^{1}\right\vert \right\} ; \\
\text{\textsc{S}}\left( J^{\ast }\right) &\equiv &\left( I\setminus J^{\ast
}\right) \setminus E\left( J^{\ast }\right) \ .
\end{eqnarray*}%
Then it suffices to show both 
\begin{eqnarray*}
A &\equiv &\sum_{J\in \mathcal{M}_{\limfunc{deep}}}\left( \frac{\mathrm{P}%
^{\alpha }\left( J,\mathbf{1}_{\text{\textsc{E}}\left( J^{\ast }\right)
}\sigma \right) }{\left\vert J\right\vert ^{\frac{1}{n}}}\right)
^{2}\left\Vert \mathsf{P}_{J}^{\limfunc{subgood},\omega }\mathbf{x}%
\right\Vert _{L^{2}\left( \omega \right) }^{2}\leq \left( \mathfrak{T}_{%
\mathbf{R}^{\alpha ,n}}^{2}+\mathcal{A}_{2}^{\alpha }\right) \left\vert
I\right\vert _{\sigma }\ , \\
B &\equiv &\sum_{J\in \mathcal{M}_{\limfunc{deep}}}\left( \frac{\mathrm{P}%
^{\alpha }\left( J,\mathbf{1}_{\text{\textsc{S}}\left( J^{\ast }\right)
}\sigma \right) }{\left\vert J\right\vert ^{\frac{1}{n}}}\right)
^{2}\left\Vert \mathsf{P}_{J}^{\limfunc{subgood},\omega }\mathbf{x}%
\right\Vert _{L^{2}\left( \omega \right) }^{2}\leq \left( \mathfrak{T}_{%
\mathbf{R}^{\alpha ,n}}^{2}+\mathcal{A}_{2}^{\alpha }\right) \left\vert
I\right\vert _{\sigma }\ .
\end{eqnarray*}%
Term $A$ is estimated in analogy with the Hilbert transform estimate (\ref%
{Hilbert}), while term $B$ is estimated by summing Poisson tails. Both
estimates rely heavily on the one-dimensional nature of $\omega $.

For $\left( x^{1},0^{\prime }\right) ,\left( z^{1},0^{\prime }\right) \in J$
in term $A$ we claim the following `strong reversal' of energy,%
\begin{eqnarray}
&&\left\vert \frac{R_{1}^{\alpha ,n}\mathbf{1}_{\text{\textsc{E}}\left(
J^{\ast }\right) }\sigma \left( x^{1},0^{\prime }\right) -R_{1}^{\alpha ,n}%
\mathbf{1}_{\text{\textsc{E}}\left( J^{\ast }\right) }\sigma \left(
z^{1},0^{\prime }\right) }{x^{1}-z^{1}}\right\vert  \label{diff quotient} \\
&=&\left\vert \int_{\text{\textsc{E}}\left( J^{\ast }\right) }\left\{ \frac{%
K_{1}^{\alpha ,n}\left( \left( x^{1},0^{\prime }\right) ,y\right)
-K_{1}^{\alpha ,n}\left( \left( z^{1},0^{\prime }\right) ,y\right) }{%
x^{1}-z^{1}}\right\} d\sigma \left( y\right) \right\vert  \notag \\
&=&\left\vert \int_{\text{\textsc{E}}\left( J^{\ast }\right) }\left\{ \frac{%
\frac{x^{1}-y^{1}}{\left( \left\vert x^{1}-y^{1}\right\vert ^{2}+\left\vert
y^{\prime }\right\vert ^{2}\right) ^{\frac{n+1-\alpha }{2}}}-\frac{%
z^{1}-y^{1}}{\left( \left\vert z^{1}-y^{1}\right\vert ^{2}+\left\vert
y^{\prime }\right\vert ^{2}\right) ^{\frac{n+1-\alpha }{2}}}}{x^{1}-z^{1}}%
\right\} d\sigma \left( y\right) \right\vert  \notag \\
&\approx &\frac{\mathrm{P}^{\alpha }\left( J,\mathbf{1}_{\text{\textsc{E}}%
\left( J^{\ast }\right) }\sigma \right) }{\left\vert J\right\vert ^{\frac{1}{%
n}}}\ .  \notag
\end{eqnarray}%
Indeed, if we set $a=\left\vert y^{\prime }\right\vert $ and $s=x^{1}-y^{1}$
and $t=z^{1}-y^{1}$, then the term in braces in (\ref{diff quotient}) is%
\begin{eqnarray*}
&&\frac{\frac{x^{1}-y^{1}}{\left( \left\vert x^{1}-y^{1}\right\vert
^{2}+\left\vert y^{\prime }\right\vert ^{2}\right) ^{\frac{n+1-\alpha }{2}}}-%
\frac{z^{1}-y^{1}}{\left( \left\vert z^{1}-y^{1}\right\vert ^{2}+\left\vert
y^{\prime }\right\vert ^{2}\right) ^{\frac{n+1-\alpha }{2}}}}{x^{1}-z^{1}} \\
&=&\frac{\frac{s}{\left( s^{2}+a^{2}\right) ^{\frac{n+1-\alpha }{2}}}-\frac{t%
}{\left( t^{2}+a^{2}\right) ^{\frac{n+1-\alpha }{2}}}}{s-t}=\frac{\varphi
\left( s\right) -\varphi \left( t\right) }{s-t}\ ,
\end{eqnarray*}%
where $\varphi \left( t\right) =t\left( t^{2}+a^{2}\right) ^{-\frac{%
n+1-\alpha }{2}}$. Now the derivative of $\varphi \left( t\right) $ is%
\begin{eqnarray*}
\frac{d}{dt}\varphi \left( t\right) &=&\left( t^{2}+a^{2}\right) ^{-\frac{%
n+1-\alpha }{2}}-\frac{n+1-\alpha }{2}\left( t^{2}+a^{2}\right) ^{-\frac{%
n+1-\alpha }{2}-1}2t^{2} \\
&=&\left( t^{2}+a^{2}\right) ^{-\frac{n+1-\alpha }{2}-1}\left\{ \left(
t^{2}+a^{2}\right) -\left( n+1-\alpha \right) t^{2}\right\} \\
&=&\left( t^{2}+a^{2}\right) ^{-\frac{n+1-\alpha }{2}-1}\left\{ a^{2}-\left(
n-\alpha \right) t^{2}\right\} ,
\end{eqnarray*}%
and since $\left\vert t\right\vert \geq \gamma \left\vert J\right\vert ^{%
\frac{1}{n}}\geq \gamma a$, we have $\left( n-\alpha \right) t^{2}\geq
\left( n-\alpha \right) \gamma ^{2}a^{2}\geq 2a^{2}$ provided we choose 
\begin{equation}
\gamma \geq \sqrt{\frac{2}{n-\alpha }}.  \label{provided'}
\end{equation}%
Thus if (\ref{provided'}) holds we get%
\begin{equation*}
-\frac{d}{dt}\varphi \left( t\right) \approx t^{2}\left( t^{2}+a^{2}\right)
^{-\frac{n+1-\alpha }{2}-1}.
\end{equation*}%
Finally, since $\left\vert s-t\right\vert \leq a\leq \frac{1}{\gamma }%
\left\vert t\right\vert \ll \left\vert t\right\vert $, the derivative $\frac{%
d\varphi }{dt}$ is essentially constant on the small interval $\left(
s,t\right) $, and we can apply the tangent line approximation to $\varphi
\left( t\right) $ to obtain $\varphi \left( s\right) -\varphi \left(
t\right) \approx \frac{d\varphi }{dt}\left( t\right) \left( s-t\right) $,
and conclude that for $\left( x^{1},0^{\prime }\right) ,\left(
z^{1},0^{\prime }\right) \in J$, 
\begin{eqnarray*}
&&\left\vert \int_{\text{\textsc{E}}\left( J^{\ast }\right) }\left\{ \frac{%
\frac{x^{1}-y^{1}}{\left( \left\vert x^{1}-y^{1}\right\vert ^{2}+\left\vert
y^{\prime }\right\vert ^{2}\right) ^{\frac{n+1-\alpha }{2}}}-\frac{%
z^{1}-y^{1}}{\left( \left\vert z^{1}-y^{1}\right\vert ^{2}+\left\vert
y^{\prime }\right\vert ^{2}\right) ^{\frac{n+1-\alpha }{2}}}}{x^{1}-z^{1}}%
\right\} d\sigma \left( y\right) \right\vert \\
&\approx &\int_{\text{\textsc{E}}\left( J^{\ast }\right) }\frac{\left\vert
z^{1}-y^{1}\right\vert ^{2}}{\left( \left\vert z^{1}-y^{1}\right\vert
^{2}+\left\vert y^{\prime }\right\vert ^{2}\right) ^{\frac{n+1-\alpha }{2}+1}%
}d\sigma \left( y\right) \approx \frac{\mathrm{P}^{\alpha }\left( J,\mathbf{1%
}_{\text{\textsc{E}}\left( J^{\ast }\right) }\sigma \right) }{\left\vert
J\right\vert ^{\frac{1}{n}}}\ ,
\end{eqnarray*}%
which proves (\ref{diff quotient}).

Thus we have%
\begin{eqnarray*}
&&\sum_{J\in \mathcal{M}_{\limfunc{deep}}}\left( \frac{\mathrm{P}^{\alpha
}\left( J,\mathbf{1}_{\text{\textsc{E}}\left( J^{\ast }\right) }\sigma
\right) }{\left\vert J\right\vert ^{\frac{1}{n}}}\right) ^{2}\int_{J\cap
L}\left\vert x^{1}-\mathbb{E}_{J}^{\omega }x^{1}\right\vert ^{2}d\omega
\left( y\right) \\
&=&\frac{1}{2}\sum_{J\in \mathcal{M}_{\limfunc{deep}}}\left( \frac{\mathrm{P}%
^{\alpha }\left( J,\mathbf{1}_{\text{\textsc{E}}\left( J^{\ast }\right)
}\sigma \right) }{\left\vert J\right\vert ^{\frac{1}{n}}}\right) ^{2}\frac{1%
}{\left\vert J\cap L\right\vert _{\omega }}\int_{J\cap L}\int_{J\cap
L}\left( x^{1}-z^{1}\right) ^{2}d\omega \left( x\right) d\omega \left(
z\right) \\
&\approx &\sum_{J\in \mathcal{M}_{\limfunc{deep}}}\frac{1}{\left\vert
J\right\vert _{\omega }}\int_{J\cap L}\int_{J\cap L}\left\{ R_{1}^{\alpha ,n}%
\mathbf{1}_{\text{\textsc{E}}\left( J^{\ast }\right) }\sigma \left(
x^{1},0^{\prime }\right) -R_{1}^{\alpha ,n}\mathbf{1}_{\text{\textsc{E}}%
\left( J^{\ast }\right) }\sigma \left( z^{1},0^{\prime }\right) \right\}
^{2}d\omega \left( x\right) d\omega \left( z\right) \\
&\lesssim &\sum_{J\in \mathcal{M}_{\limfunc{deep}}}\frac{1}{\left\vert
J\right\vert _{\omega }}\int_{J\cap L}\int_{J\cap L}\left\{ R_{1}^{\alpha ,n}%
\mathbf{1}_{I}\sigma \left( x^{1},0^{\prime }\right) -R_{1}^{\alpha ,n}%
\mathbf{1}_{I}\sigma \left( z^{1},0^{\prime }\right) \right\} ^{2}d\omega
\left( x\right) d\omega \left( z\right) \\
&&+\sum_{J\in \mathcal{M}_{\limfunc{deep}}}\frac{1}{\left\vert J\right\vert
_{\omega }}\int_{J\cap L}\int_{J\cap L}\left\{ R_{1}^{\alpha ,n}\mathbf{1}%
_{J^{\ast }}\sigma \left( x^{1},0^{\prime }\right) -R_{1}^{\alpha ,n}\mathbf{%
1}_{J^{\ast }}\sigma \left( z^{1},0^{\prime }\right) \right\} ^{2}d\omega
\left( x\right) d\omega \left( z\right) \\
&&+\sum_{J\in \mathcal{M}_{\limfunc{deep}}}\frac{1}{\left\vert J\right\vert
_{\omega }}\int_{J\cap L}\int_{J\cap L}\left\{ R_{1}^{\alpha ,n}\mathbf{1}_{%
\text{\textsc{S}}\left( J^{\ast }\right) }\sigma \left( x^{1},0^{\prime
}\right) -R_{1}^{\alpha ,n}\mathbf{1}_{\text{\textsc{S}}\left( J^{\ast
}\right) }\sigma \left( z^{1},0^{\prime }\right) \right\} ^{2}d\omega \left(
x\right) d\omega \left( z\right) \\
&\equiv &A_{1}+A_{2}+A_{3},
\end{eqnarray*}%
since $I=J^{\ast }\dot{\cup}\left( I\setminus J^{\ast }\right) =J^{\ast }%
\dot{\cup}$\textsc{E}$\left( J^{\ast }\right) \dot{\cup}$\textsc{S}$\left(
J^{\ast }\right) $. Now we can discard the difference in term $A_{1}$ by
writing 
\begin{equation*}
\left\vert R_{1}^{\alpha ,n}\mathbf{1}_{I}\sigma \left( x^{1},0^{\prime
}\right) -R_{1}^{\alpha ,n}\mathbf{1}_{I}\sigma \left( z^{1},0^{\prime
}\right) \right\vert \leq \left\vert R_{1}^{\alpha ,n}\mathbf{1}_{I}\sigma
\left( x^{1},0^{\prime }\right) \right\vert +\left\vert R_{1}^{\alpha ,n}%
\mathbf{1}_{I}\sigma \left( z^{1},0^{\prime }\right) \right\vert
\end{equation*}%
to obtain from pairwise disjointedness of $J\in \mathcal{M}_{\limfunc{deep}}$%
, 
\begin{equation*}
A_{1}\lesssim \sum_{J\in \mathcal{M}_{\limfunc{deep}}}\int_{J\cap
L}\left\vert R_{1}^{\alpha ,n}\mathbf{1}_{I}\sigma \left( x^{1},0^{\prime
}\right) \right\vert ^{2}d\omega \left( x\right) \leq \int_{I}\left\vert
R_{1}^{\alpha ,n}\mathbf{1}_{I}\sigma \right\vert ^{2}d\omega \leq \mathfrak{%
T}_{R_{1}^{\alpha ,n}}^{2}\left\vert I\right\vert _{\sigma }\ ,
\end{equation*}%
and similarly we can discard the difference in term $A_{2}$, and use the
bounded overlap property (\ref{bounded overlap}), to obtain%
\begin{eqnarray*}
A_{2} &\lesssim &\sum_{J\in \mathcal{M}_{\limfunc{deep}}}\int_{J\cap
L}\left\vert R_{1}^{\alpha ,n}\mathbf{1}_{J^{\ast }}\sigma \left(
x^{1},0^{\prime }\right) \right\vert ^{2}d\omega \left( x\right) \leq
\sum_{J\in \mathcal{M}_{\limfunc{deep}}}\mathfrak{T}_{R_{1}^{\alpha
,n}}^{2}\left\vert J^{\ast }\right\vert _{\sigma } \\
&=&\mathfrak{T}_{R_{1}^{\alpha ,n}}^{2}\sum_{r=1}^{\infty }\sum_{J\in 
\mathcal{M}_{\limfunc{deep}}\left( I_{r}\right) }\left\vert J^{\ast
}\right\vert _{\sigma }\leq \mathfrak{T}_{R_{1}^{\alpha
,n}}^{2}\sum_{r=1}^{\infty }\beta \left\vert I_{r}\right\vert _{\sigma }\leq
\beta \mathfrak{T}_{R_{1}^{\alpha ,n}}^{2}\left\vert I\right\vert _{\sigma
}\ .
\end{eqnarray*}

\begin{remark}
\label{fails}The above estimate fails for the \emph{nearby} cubes $J$ in $%
I_{r}$, and so it is important to use the definition of the energy condition
as in Definition \ref{energy condition}\ above.
\end{remark}

This leaves us to consider the term%
\begin{eqnarray*}
A_{3} &=&\sum_{J\in \mathcal{M}_{\limfunc{deep}}}\frac{1}{\left\vert
J\right\vert _{\omega }}\int_{J\cap L}\int_{J\cap L}\left\{ R_{1}^{\alpha ,n}%
\mathbf{1}_{\text{\textsc{S}}\left( J^{\ast }\right) }\sigma \left(
x^{1},0^{\prime }\right) -R_{1}^{\alpha ,n}\mathbf{1}_{\text{\textsc{S}}%
\left( J^{\ast }\right) }\sigma \left( z^{1},0^{\prime }\right) \right\}
^{2}d\omega \left( x\right) d\omega \left( z\right) \\
&=&2\sum_{J\in \mathcal{M}_{\limfunc{deep}}}\int_{J\cap L}\left\{
R_{1}^{\alpha ,n}\mathbf{1}_{\text{\textsc{S}}\left( J^{\ast }\right)
}\sigma \left( x^{1},0^{\prime }\right) -\mathbb{E}_{J\cap L}^{\omega }\left[
R_{1}^{\alpha ,n}\mathbf{1}_{\text{\textsc{S}}\left( J^{\ast }\right)
}\sigma \left( z^{1},0^{\prime }\right) \right] \right\} ^{2}d\omega \left(
x\right) ,
\end{eqnarray*}%
in which we do \emph{not} discard the difference. However, because the
average is subtracted off, we can apply the Energy Lemma \ref{ener} to each
term in this sum to dominate it by,%
\begin{equation}
B=\sum_{J\in \mathcal{M}_{\limfunc{deep}}}\left( \frac{\mathrm{P}^{\alpha
}\left( J,\mathbf{1}_{\text{\textsc{S}}\left( J^{\ast }\right) }\sigma
\right) }{\left\vert J\right\vert ^{\frac{1}{n}}}\right) ^{2}\left\Vert 
\mathsf{P}_{\mathcal{J}}^{\omega }\mathbf{x}\right\Vert _{L^{2}\left( \omega
\right) }^{2}\ .  \label{B1 + B2}
\end{equation}

To estimate $B$, we first assume that $n-1\leq \alpha <n$ so that $\mathrm{P}%
^{\alpha }\left( J,\mathbf{1}_{\text{\textsc{S}}\left( J^{\ast }\right)
}\sigma \right) \leq \mathcal{P}^{\alpha }\left( J,\mathbf{1}_{\text{\textsc{%
S}}\left( J^{\ast }\right) }\sigma \right) $, and then use $\left\Vert 
\mathsf{P}_{\mathcal{J}}^{\omega }\mathbf{x}\right\Vert _{L^{2}\left( \omega
\right) }^{2}\leq \left\vert J\right\vert ^{\frac{2}{n}}\left\vert
J\right\vert _{\omega }$ and apply the $\mathcal{A}_{2}^{\alpha }$ condition
to obtain the following `pivotal reversal' of energy,%
\begin{eqnarray*}
B &\leq &\sum_{J\in \mathcal{M}_{\limfunc{deep}}}\mathrm{P}^{\alpha }\left(
J,\mathbf{1}_{\text{\textsc{S}}\left( J^{\ast }\right) }\sigma \right)
^{2}\left\vert J\right\vert _{\omega }\leq \sum_{J\in \mathcal{M}_{\limfunc{%
deep}}}\mathrm{P}^{\alpha }\left( J,\mathbf{1}_{\text{\textsc{S}}\left(
J^{\ast }\right) }\sigma \right) \left\{ \mathcal{P}^{\alpha }\left( J,%
\mathbf{1}_{\text{\textsc{S}}\left( J^{\ast }\right) }\sigma \right)
\left\vert J\right\vert _{\omega }\right\} \\
&\leq &\mathcal{A}_{2}^{\alpha }\sum_{J\in \mathcal{M}_{\limfunc{deep}}}%
\mathrm{P}^{\alpha }\left( J,\mathbf{1}_{\text{\textsc{S}}\left( J^{\ast
}\right) }\sigma \right) \left\vert J\right\vert ^{1-\frac{\alpha }{n}}=%
\mathcal{A}_{2}^{\alpha }\sum_{J\in \mathcal{M}_{\limfunc{deep}}}\int_{\text{%
\textsc{S}}\left( J^{\ast }\right) }\frac{\left\vert J\right\vert ^{\frac{1}{%
n}}\left\vert J\right\vert ^{1-\frac{\alpha }{n}}}{\left( \left\vert
J\right\vert ^{\frac{1}{n}}+\left\vert y-c_{J}\right\vert \right)
^{n+1-\alpha }}d\sigma \left( y\right) \\
&=&\mathcal{A}_{2}^{\alpha }\sum_{J\in \mathcal{M}_{\limfunc{deep}}}\int_{%
\text{\textsc{S}}\left( J^{\ast }\right) }\left( \frac{\left\vert
J\right\vert ^{\frac{1}{n}}}{\left\vert J\right\vert ^{\frac{1}{n}%
}+\left\vert y-c_{J}\right\vert }\right) ^{n+1-\alpha }d\sigma \left(
y\right) \\
&=&\mathcal{A}_{2}^{\alpha }\int_{I}\left\{ \sum_{J\in \mathcal{M}_{\limfunc{%
deep}}}\left( \frac{\left\vert J\right\vert ^{\frac{1}{n}}}{\left\vert
J\right\vert ^{\frac{1}{n}}+\left\vert y-c_{J}\right\vert }\right)
^{n+1-\alpha }\mathbf{1}_{\text{\textsc{S}}\left( J^{\ast }\right) }\left(
y\right) \right\} d\sigma \left( y\right) \\
&\equiv &\mathcal{A}_{2}^{\alpha }\int_{I}F\left( y\right) d\sigma \left(
y\right) .
\end{eqnarray*}%
At this point we claim that $F\left( y\right) \leq C$ with a constant $C$
independent of the decomposition $\mathcal{M}_{\limfunc{deep}}=\overset{%
\cdot }{\dbigcup\limits_{r\geq 1}}\mathcal{M}_{\limfunc{deep}}\left(
I_{r}\right) $. Indeed, if $y$ is fixed, then the only cubes $J\in \mathcal{M%
}_{\limfunc{deep}}$ for which $y\in $\textsc{S}$\left( J^{\ast }\right) $
are those $J$ satisfying%
\begin{equation*}
J\cap \limfunc{Sh}\left( y;\gamma \right) \neq \emptyset ,
\end{equation*}%
where $\limfunc{Sh}\left( y;\gamma \right) $ is the Carleson shadow of the
point $y$ onto the $x_{1}$-axis $L$ with sides of slope $\frac{1}{\gamma }$,
i.e. $\limfunc{Sh}\left( y;\gamma \right) $ is interval on $L$ with length $%
2\gamma \limfunc{dist}\left( y,L\right) $ and center equal to the point on $%
L $ that is closest to $y$. Now there can be at most two cubes $J$ whose
side length exceeds $2\gamma \limfunc{dist}\left( y,L\right) $, and for
these cubes we simply use $\frac{\left\vert J\right\vert ^{\frac{1}{n}}}{%
\left\vert J\right\vert ^{\frac{1}{n}}+\left\vert y-c_{J}\right\vert }\leq 1$%
. As for the remaining cubes $J$, they are all contained inside the triple $3%
\limfunc{Sh}\left( y;\gamma \right) $ of the shadow, and the distance $%
\left\vert y-c_{J}\right\vert $ is essentially $\limfunc{dist}\left(
y,L\right) $ (up to a factor of $\gamma $) for all of these cubes. Thus we
have the estimate%
\begin{eqnarray*}
\sum_{\substack{ J\in \mathcal{M}_{\limfunc{deep}}  \\ J\subset 3\limfunc{Sh}%
\left( y;\gamma \right) }}\left( \frac{\left\vert J\right\vert ^{\frac{1}{n}}%
}{\left\vert J\right\vert ^{\frac{1}{n}}+\left\vert y-c_{J}\right\vert }%
\right) ^{n+1-\alpha } &\lesssim &\sum_{\substack{ J\in \mathcal{M}_{%
\limfunc{deep}}  \\ J\subset 3\limfunc{Sh}\left( y;\gamma \right) }}\left( 
\frac{\left\vert J\right\vert ^{\frac{1}{n}}}{\limfunc{dist}\left(
y,L\right) }\right) ^{n+1-\alpha } \\
&\lesssim &\frac{1}{\limfunc{dist}\left( y,L\right) ^{n+1-\alpha }}\sum 
_{\substack{ J\in \mathcal{M}_{\limfunc{deep}}  \\ J\subset 3\limfunc{Sh}%
\left( y;\gamma \right) }}\left\vert J\cap L\right\vert ^{n+1-\alpha } \\
&\lesssim &\frac{1}{\limfunc{dist}\left( y,L\right) ^{n+1-\alpha }}\sum 
_{\substack{ J\in \mathcal{M}_{\limfunc{deep}}  \\ J\subset 3\limfunc{Sh}%
\left( y;\gamma \right) }}\limfunc{dist}\left( y,L\right) ^{n-\alpha }\
\left\vert J\cap L\right\vert \\
&\lesssim &\frac{\limfunc{dist}\left( y,L\right) ^{n-\alpha }}{\limfunc{dist}%
\left( y,L\right) ^{n+1-\alpha }}\left\vert 3\limfunc{Sh}\left( y;\gamma
\right) \right\vert \lesssim 1,
\end{eqnarray*}%
because the intervals $\left\{ J\cap L\right\} _{\substack{ J\in \mathcal{M}%
_{\limfunc{deep}}  \\ J\subset 3\limfunc{Sh}\left( y;\gamma \right) }}$ are
pairwise disjoint in $3\limfunc{Sh}\left( y;\gamma \right) $, and $%
\left\vert J\cap L\right\vert ^{\frac{1}{n}}$ is the length of $J\cap L$,
and since $n+1-\alpha >1=\dim L$. It is here that the one-dimensional nature
of $\omega $ delivers the boundedness of this sum of Poisson tails. Thus we
have%
\begin{equation*}
B\leq \mathcal{A}_{2}^{\alpha }\int_{I}F\left( y\right) d\sigma \left(
y\right) \leq C\mathcal{A}_{2}^{\alpha }\left\vert I\right\vert _{\sigma }\ ,
\end{equation*}%
which is the desired estimate in the case that $n-1\leq \alpha <n$.

Now we suppose that $0\leq \alpha <n-1$ and use Cauchy-Schwarz to obtain%
\begin{eqnarray*}
\mathrm{P}^{\alpha }\left( J,\mathbf{1}_{\text{\textsc{S}}\left( J^{\ast
}\right) }\sigma \right) &=&\int_{\text{\textsc{S}}\left( J^{\ast }\right) }%
\frac{\left\vert J\right\vert ^{\frac{1}{n}}}{\left( \left\vert J\right\vert
^{\frac{1}{n}}+\left\vert y-c_{J}\right\vert \right) ^{n+1-\alpha }}d\sigma
\left( y\right) \\
&\leq &\left\{ \int_{\text{\textsc{S}}\left( J^{\ast }\right) }\frac{%
\left\vert J\right\vert ^{\frac{1}{n}}}{\left( \left\vert J\right\vert ^{%
\frac{1}{n}}+\left\vert y-c_{J}\right\vert \right) ^{n+1-\alpha }}\left( 
\frac{\left\vert J\right\vert ^{\frac{1}{n}}}{\left\vert J\right\vert ^{%
\frac{1}{n}}+\left\vert y-c_{J}\right\vert }\right) ^{n-1-\alpha }d\sigma
\left( y\right) \right\} ^{\frac{1}{2}} \\
&&\times \left\{ \int_{\text{\textsc{S}}\left( J^{\ast }\right) }\frac{%
\left\vert J\right\vert ^{\frac{1}{n}}}{\left( \left\vert J\right\vert ^{%
\frac{1}{n}}+\left\vert y-c_{J}\right\vert \right) ^{n+1-\alpha }}\left( 
\frac{\left\vert J\right\vert ^{\frac{1}{n}}}{\left\vert J\right\vert ^{%
\frac{1}{n}}+\left\vert y-c_{J}\right\vert }\right) ^{\alpha +1-n}d\sigma
\left( y\right) \right\} ^{\frac{1}{2}} \\
&=&\mathcal{P}^{\alpha }\left( J,\mathbf{1}_{\text{\textsc{S}}\left( J^{\ast
}\right) }\sigma \right) ^{\frac{1}{2}} \\
&&\times \left\{ \int_{\text{\textsc{S}}\left( J^{\ast }\right) }\frac{%
\left( \left\vert J\right\vert ^{\frac{1}{n}}\right) ^{\alpha +2-n}}{\left(
\left\vert J\right\vert ^{\frac{1}{n}}+\left\vert y-c_{J}\right\vert \right)
^{2}}d\sigma \left( y\right) \right\} ^{\frac{1}{2}}.
\end{eqnarray*}%
Then arguing as above we have%
\begin{eqnarray*}
B &\leq &\sum_{J\in \mathcal{M}_{\limfunc{deep}}}\mathrm{P}^{\alpha }\left(
J^{\ast },\mathbf{1}_{\text{\textsc{S}}\left( J^{\ast }\right) }\sigma
\right) ^{2}\left\vert J\right\vert _{\omega } \\
&\leq &\sum_{J\in \mathcal{M}_{\limfunc{deep}}}\left\{ \mathcal{P}^{\alpha
}\left( J^{\ast },\mathbf{1}_{\text{\textsc{S}}\left( J^{\ast }\right)
}\sigma \right) \left\vert J\right\vert _{\omega }\right\} \int_{\text{%
\textsc{S}}\left( J^{\ast }\right) }\frac{\left( \left\vert J\right\vert ^{%
\frac{1}{n}}\right) ^{\alpha +2-n}}{\left( \left\vert J\right\vert ^{\frac{1%
}{n}}+\left\vert y-c_{J}\right\vert \right) ^{2}}d\sigma \left( y\right) \\
&\leq &\mathcal{A}_{2}^{\alpha }\sum_{J\in \mathcal{M}_{\limfunc{deep}%
}}\left\vert J\right\vert ^{1-\frac{\alpha }{n}}\int_{\text{\textsc{S}}%
\left( J^{\ast }\right) }\frac{\left( \left\vert J\right\vert ^{\frac{1}{n}%
}\right) ^{\alpha +2-n}}{\left( \left\vert J\right\vert ^{\frac{1}{n}%
}+\left\vert y-c_{J}\right\vert \right) ^{2}}d\sigma \left( y\right) =%
\mathcal{A}_{2}^{\alpha }\sum_{J\in \mathcal{M}_{\limfunc{deep}}}\int_{\text{%
\textsc{S}}\left( J^{\ast }\right) }\frac{\left\vert J\right\vert ^{\frac{2}{%
n}}}{\left( \left\vert J\right\vert ^{\frac{1}{n}}+\left\vert
y-c_{J}\right\vert \right) ^{2}}d\sigma \left( y\right) \\
&=&\mathcal{A}_{2}^{\alpha }\int_{I}\left\{ \sum_{J\in \mathcal{M}_{\limfunc{%
deep}}}\left( \frac{\left\vert J\right\vert ^{\frac{1}{n}}}{\left\vert
J\right\vert ^{\frac{1}{n}}+\left\vert y-c_{J}\right\vert }\right) ^{2}%
\mathbf{1}_{\text{\textsc{S}}\left( J^{\ast }\right) }\left( y\right)
\right\} d\sigma \left( y\right) \equiv \mathcal{A}_{2}^{\alpha
}\int_{I}F\left( y\right) d\sigma \left( y\right) ,
\end{eqnarray*}%
and again $F\left( y\right) \leq C$ because 
\begin{eqnarray*}
\sum_{\substack{ J\in \mathcal{M}_{\limfunc{deep}}  \\ J\subset 3\limfunc{Sh}%
\left( y;\gamma \right) }}\left( \frac{\left\vert J\right\vert ^{\frac{1}{n}}%
}{\left\vert J\right\vert ^{\frac{1}{n}}+\left\vert y-c_{J}\right\vert }%
\right) ^{2} &\lesssim &\frac{1}{\limfunc{dist}\left( y,L\right) ^{2}}\sum 
_{\substack{ J\in \mathcal{M}_{\limfunc{deep}}  \\ J\subset 3\limfunc{Sh}%
\left( y;\gamma \right) }}\left\vert J\right\vert ^{\frac{2}{n}} \\
&\lesssim &\frac{1}{\limfunc{dist}\left( y,L\right) ^{2}}\sum_{\substack{ %
J\in \mathcal{M}_{\limfunc{deep}}  \\ J\subset 3\limfunc{Sh}\left( y;\gamma
\right) }}\left\vert J\cap L\right\vert ^{2} \\
&\lesssim &\frac{1}{\limfunc{dist}\left( y,L\right) ^{2}}\sum_{\substack{ %
J\in \mathcal{M}_{\limfunc{deep}}  \\ J\subset 3\limfunc{Sh}\left( y;\gamma
\right) }}\limfunc{dist}\left( y,L\right) \ \left\vert J\cap L\right\vert \\
&\lesssim &\frac{\limfunc{dist}\left( y,L\right) }{\limfunc{dist}\left(
y,L\right) ^{2}}\left\vert 3\limfunc{Sh}\left( y;\gamma \right) \right\vert
\lesssim 1.
\end{eqnarray*}%
Thus we again have%
\begin{equation*}
B\leq \mathcal{A}_{2}^{\alpha }\int_{I}F\left( y\right) d\sigma \left(
y\right) \leq C\mathcal{A}_{2}^{\alpha }\left\vert I\right\vert _{\sigma }\ ,
\end{equation*}%
and this completes the proof of necessity of the energy conditions when one
of the measures is supported on a line.

\end{document}